\documentclass[11pt]{article}
\usepackage{amsfonts}

\usepackage{graphics}
\usepackage{cite}
\usepackage{latexsym}
\usepackage{amsmath}
\usepackage{amssymb}
\usepackage[dvips]{epsfig}
\usepackage{amscd}

\newtheorem{theorem}{Theorem}[section]
\newtheorem{remark}{Remark}[section]

\newtheorem{lemma}[theorem]{Lemma}

\newtheorem{proposition}[theorem]{Proposition}

\newcommand{\n}{\rho}

\newcommand{\ti}{\tilde}

\def\pf{{\it Proof.}  }
\renewcommand{\div}{ {\rm div }  }

\newcommand{\na}{\nabla }

\newcommand{\pa}{\partial}

\newcommand{\bt}{\begin{theorem}}
\newcommand{\bl}{\begin{lemma}}
\newcommand{\el}{\end{lemma}}
\newcommand{\et}{\end{theorem}}

\newcommand{\ve}{\varepsilon}
\newcommand{\la}{\label}
\newcommand{\si}{\sigma}

\newcommand{\bn}{\begin{eqnarray}}
\newcommand{\en}{\end{eqnarray}}
\newcommand{\bnn}{\begin{eqnarray*}}
\newcommand{\enn}{\end{eqnarray*}}

\newcommand{\bnnn}{\begin{eqnarray*}}
\newcommand{\ennn}{\end{eqnarray*}}
\newcommand{\ben}{\begin{enumerate}}
\newcommand{\een}{\end{enumerate}}

\newcommand{\ba}{\begin{aligned}}
\newcommand{\ea}{\end{aligned}}
\newcommand{\be}{\begin{equation}}
\newcommand{\ee}{\end{equation}}

\def\norm[#1]#2{\|#2\|_{#1}}

\def\rrr{\mathbb{R}^3}

\def\r{\mathbb{R}}

\newcommand{\mn}{\mu(\rho)}
\newcommand{\ban}{\bar\rho}

\newcommand{\xmu}{\underline{\mu}}
\newcommand{\smu}{\bar\mu}

\newcommand{\xl}{\left}
\newcommand{\xr}{\right}

\makeatletter
\@addtoreset{equation}{section}
\makeatother

\makeatletter
\@addtoreset{equation}{section}
\makeatother

\title{On  the  Cauchy Problem of 3D   Nonhomogeneous Navier-Stokes Equations with   Density-Dependent Viscosity and Vacuum
}

\author{Cheng He $^{\rm a},$ Jing Li $^{\rm b,c}$, Boqiang L\"u $^{\rm d}$ \thanks{
  The research of \textsc{J. Li} was
partially supported  by the National Center for Mathematics and Interdisciplinary Sciences, CAS, and  NNSFC Grant Nos.   11371348, 11688101,  and 11525106.    The research of \textsc{B. L\"u} was
partially   supported by NNSFC (No. 11601218) and Natural Science Foundation of Jiangxi Province (No. 20161BAB211002).
 Email: hecheng@nsfc.gov.cn (C. He),  ajingli@gmail.com (J. Li), lvbq86@163.com(B. L\"u).}\\
{\normalsize a. Division of Mathematics, Department of Mathematical  \&  Physical Sciences,}\\{\normalsize  National Natural Science Foundation of China,} \\{\normalsize  Beijing 100085, P. R. China} \\
  {\normalsize b. NCMIS, HLM, Institute of Applied Mathematics,}\\
  Academy of Mathematics and Systems Science,  \\
{\normalsize  Chinese Academy of Sciences,    Beijing 100190,
P. R. China}  \\ c. School of Mathematical Sciences, University of Chinese Academy of Sciences,\\  Beijing 100049, P. R. China\\
  {\normalsize   d.  College of Mathematics and Information Science, } \\
   {\normalsize  Nanchang Hangkong University,} \\ {\normalsize   Nanchang 330063, P. R. China} }

\date{ }

\begin{document}
\maketitle

\begin{abstract}
We consider the global existence and large-time asymptotic behavior  of strong solutions to the Cauchy problem of the three-dimensional   nonhomogeneous incompressible Navier-Stokes equations with density-dependent viscosity and vacuum.
 %($\|u_0\|_{\dot H^\beta}$ with $\beta\in(1/2, 1]$),
We   establish some key a priori   exponential   decay-in-time rates of  the strong solutions.  Then after using these estimates, we also obtain the global existence   of strong solutions     in the whole three-dimensional space, provided  that the  initial velocity    is suitably small in the $
\dot H^\beta$-norm for some $\beta\in(1/2,1].$     Note that  this result is proved without any smallness conditions on the initial density. Moreover, the density can contain vacuum states and even have compact support initially.
%In addition, when the viscosity is a positive constant, we also prove some new time-independent estimates and  the large time  exponential  decay rates for the global strong solution under the smallness assumption on $\|u_0\|_{\dot H^{1/2}}$.
\end{abstract}


Math Subject Classification: 35Q35; 76N10.

\section{Introduction}

The   nonhomogeneous incompressible   Navier-Stokes equations (\cite{L1996}) read as follows:
\be\label{1.1}
\begin{cases}
\partial_{t}\rho+\div(\rho u)=0, \\
\partial_{t}(\rho u )+\div(\rho u \otimes u )-\div(2\mn d)  +\nabla P=0, \\
\div u=0.
\end{cases}
\ee
Here, $t\ge 0$ is time, $x\in\Omega \subset \r^3$ is the spatial coordinates, and the unknown functions $\rho=\rho(x,t)$, $u=(u^1,u^2,u^3)(x,t)$, and $P=P(x,t)$ denote the density, velocity, and pressure of the fluid, respectively; The deformation tensor is defined by
\be\la{d}
d=\frac{1}{2}\xl[\na u+(\na u)^T\xr],\ee
and the viscosity   $\mu(\n)$     satisfies  the following hypothesis:
\be\la{n3}
\mu \in C^1[0,\infty),~~   \mu(\n)>0.
\ee

We consider the Cauchy problem of \eqref{1.1} with $(\rho, u )$ vanishing at infinity  and the initial conditions:
\begin{equation}\label{1.3}
\rho(x,0)=\rho_0(x),\  \n u (x,0)=m_0(x), \ x\in\mathbb{R}^3,
\end{equation}
for given initial data $\rho_0$ and $m_0$.

%The study of the system \eqref{1.1} has a long history. Since the pioneering work of Leray \cite{L1934}, there are huge literatures on the studies of the large time existence and behavior of solutions of \eqref{1.1} in the case that the density is constant (refer to \cite{CF1988,K1984,L1969,L1996,T2001} and references therein).

There are   lots of literatures on the mathematical study of nonhomogeneous incompressible flow.
In particular,  the system \eqref{1.1} with constant viscosity  has been considered extensively.
On the one hand, in the absence of vacuum, the global existence of weak solutions and the local existence of strong ones were established in Kazhikov \cite{K1974,AK1973}. Ladyzhenskaya-Solonnikov \cite{lady} first proved the global well-posedness of strong solutions to the  initial boundary value problems in both 2D bounded domains (for large data)
and 3D ones (with initial velocity small in  suitable norms).
Recently, the global well-posedness results   with small initial data
in critical spaces were considered by many people (see \cite{dan2,zhang1,dan3,zhang3} %\cite{dan1,dan2,zhang1,zhang2,dan3,dan4,zhang3}
and the references therein).
 On the other hand, when the initial density allows to vanish,   the global existence of weak solutions is proved by Simon \cite{S1990}.   The local existence of strong solutions was obtained by Choe-Kim \cite{kim2003} (for 3D bounded and unbounded domains) and L\"u-Xu-Zhong  \cite{lvshizh1} (for 2D Cauchy problem) under some compatibility conditions. Recently, for the Cauchy problem in the whole 2D space,      L\"u-Shi-Zhong \cite{lvshizh} obtained    the global  strong solutions for large initial data. For the 3D case,   under some smallness conditions on the initial velocity,   Craig-Huang-Wang \cite{huang13}     proved the following interesting result.

 \begin{proposition}[\cite{huang13}]\label{thm3} Let $\Omega=\r^3.$ For positive constants $\bar\n $   and $\mu $,
assume that $\mn\equiv\mu$ in \eqref{1.1} and the initial data $(\rho_0, m_0)$ satisfy
\begin{equation}\label{hw2.2}
\begin{cases}
0\le \rho_{0}\le \bar\n,\  \rho_{0} \in  L^{3/2}(\r^3) \cap H^1(\r^3),\\  u_{0}\in \dot H^{1/2}(\r^3) \cap D_{0,\sigma}^1(\r^3) \cap D^{2,2}(\r^3),\ m_0=\n_0u_0,
\end{cases}\end{equation}
and the compatibility condition
\begin{equation}\label{hw2.1}
-\mu\Delta u_0+\na P_0=\n_0^{1/2}g, ~~~~~\mbox{in}~~\rrr,
\end{equation}
for some $(P_0,g)\in D^1(\rrr)\times L^2(\rrr)$. Then, there exists some positive constant $\ve$ depending only on $\bar \n$ such that  there exists a unique  global strong solution to the Cauchy problem \eqref{1.1}  \eqref{1.3}  provided $\|u_0\|_{\dot H^{1/2}}\le \mu\ve.$  Moreover,   the following large time decay rate  holds for $t\ge 1,$
\be \label{hw1}
\|\na  u(\cdot,t)\|_{L^2(\rrr)} \le \bar{C}t^{-1/2},
\ee
where $\bar C$ depends on $\ban,~\mu,$ and $\|\n_0^{1/2}u_0\|_{L^2(\rrr)}$.
\end{proposition}

When it comes to the case that the  viscosity $\mn$ depends on the density $\n$,
it is more difficult   to investigate the global well-posedness of system \eqref{1.1} due to the strong coupling between viscosity coefficient and density.
%Next,  we will list some results of the general density-dependent viscosity case.
 In fact, allowing the density to vanish initially, Lions \cite{L1996} first obtained the global weak solutions whose   uniqueness and regularity   are still open even in two spatial  dimensions. Later, Desjardins \cite{D1997} established  the global weak  solution with higher  regularity  for
2D case provided  that  the viscosity  $\mn$ is a small perturbation of a positive constant in $L^\infty$-norm. Recently,
some  progress  has been made on the well-posedness of  strong solutions to \eqref{1.1}    (see \cite{zhang5,zhang6,kim2004,HW2014,huang15,zjw,hls} and the reference therein). In particular, on the one hand,  when the initial density is strictly away from vacuum, Abidi-Zhang \cite{zhang5} obtained the global strong solutions in whole 2D space under  smallness conditions on $\|\mu(\n_0)-1\|_{L^\infty}$, and later for 3D case,  they \cite{zhang6}  obtained the global strong ones under the smallness conditions on both
 $\|u_0\|_{L^2}\|\na u_0\|_{L^2}$ and $\|\mu(\n_0)-1\|_{L^\infty}$.
On the other hand,  for the case that the initial density contains vacuum,
  Huang-Wang \cite{HW2014} obtained the global strong solutions  in 2D bounded domains when $\|\na\mu(\n_0)\|_{L^p}(p\ge2)$ is small enough;    Huang-Wang \cite{huang15} and Zhang \cite{zjw} established  the global strong solutions with small $\|\na u_0\|_{L^2}$ in 3D bounded domains. However, as pointed by  Huang-Wang \cite{huang15}, the methods used in  \cite{huang15,zjw} depend  heavily on the   boundedness of the domains and  little is known for the global well-posedness of strong solutions to the Cauchy problem  \eqref{1.1}--\eqref{1.3}    with density-dependent viscosity   and vacuum. % as far-field density.
 %As mentioned in Huang-Wang \cite{huang15},  their arguments are  only valid for  the  case that the domains are   bounded  and can not be applied  to the Cauchy problem in the whole 3D  space.

Before stating the main results, we first explain the notations and
conventions used throughout this paper. Set
$$ \int fdx\triangleq\int_{\r^3}fdx.$$ Moreover, for $1\le r\le \infty, k\ge 1, $ and $\beta>0,$ the standard homogeneous and inhomogeneous Sobolev spaces are defined as follows:
   \bnn  \begin{cases}L^r=L^r(\r^3 ),\quad
W^{k,r}  = W^{k,r}(\r^3) , \quad H^k = W^{k,2} ,\\ \|\cdot\|_{B_1\cap B_2}=\|\cdot\|_{B_1 }+\|\cdot\|_{B_2}, \mbox{ for two Banach spaces } B_1 \mbox{ and } B_2, \\
 D^{k,r}=D^{k,r}(\r^3)=\{v\in L^1_{\rm loc}(\r^3)| \na^k v\in L^r(\r^3)\},\\
D^1  =\{v\in L^6 (\r^3)| \na  v\in L^2(\r^3)\},   \\
 C_{0,\sigma}^\infty=\{f\in C_0^\infty~|~\div f=0\},\quad D_{0,\sigma}^1=\overline{C_{0,\sigma}^\infty}~\mbox{closure~in~the~norm~of}~D^{1},\\
 \dot H^\beta=\left\{f:\r^3
 \rightarrow \r\left|\|f\|^2_{\dot H^\beta}=
 \displaystyle{\int} |\xi|^{2\beta}|\hat f(\xi)|^2d\xi<\infty\right.
 \right\} ,\end{cases}\enn  where $\hat f$ is the Fourier transform
 of $f.$

%Now we define precisely what we mean by strong solutions. \begin{definition}\label{def1} If all derivatives involved in \eqref{1.1} for $(\rho,u,P)$ are regular distributions, and equations \eqref{1.1} hold almost everywhere in $\rrr\times(0,T)$, then $(\rho,u,P)$ is called a strong solution to \eqref{1.1}.\end{definition}

Our main result can be stated as follows:
\begin{theorem}\label{thm1} For constants $\bar\n>0,$ $q\in(3,\infty),$ and $\beta\in(\frac{1}{2}, 1]$,
assume that the initial data $(\rho_0, m_0)$ satisfy
\begin{equation}\label{2.2}
0\le \rho_{0}\le \bar\n,\  \rho_{0} \in L^{3/2}\cap H^1, \  \na\mu(\n_0)\in L^q  ,\  u_{0}\in \dot H^\beta\cap D_{0,\sigma}^1 ,\ m_0=\n_0u_0.
\end{equation}
% and the compatibility condition \begin{equation}\label{2.1} -\div (\mu(\n_0)\xl[\na u_0+(\na u_0)^T\xr])+\na P_0=\n_0^{1/2}g, ~~~~~\mbox{ in } \rrr, \end{equation} for some $(P_0,g)\in D^1\times L^2$.
Then for $$ \xmu\triangleq\min_{\n\in[0,\bar\n]}\mn, \quad\smu\triangleq\max_{\n\in[0,\bar\n]}\mn, \quad M\triangleq\|\na \mu(\n_0)\|_{L^q},$$  there exists some small positive constant $\ve_0$ depending only on $q, \beta, \bar\n, \xmu, \smu , \|\n_0\|_{L^{3/2}},$ and $M$ such that if
\begin{equation}\label{xx}
\|u_0\|_{\dot H^\beta}\le \ve_0,\end{equation}  the Cauchy problem \eqref{1.1}--\eqref{1.3} admits a unique global strong solution $(\rho, u, P)$ satisfying that
for any $0<\tau<T<\infty$ and $p\in [2,p_0)$  with $p_0 \triangleq\min\{6, q\},$
 \be \label{2.3} \begin{cases}
0\le \rho\in C([0,T]; L^{3/2}\cap H^1 ),
\quad \na\mn \in  C([0,T]; L^q) ,  \\
 \na u  \in L^\infty(0,T;L^2)\cap L^\infty(\tau,T;  W^{1,p_0})\cap C([\tau,T]; H^1\cap W^{1,p })  , \\
 P\in  L^\infty(\tau,T;  W^{1,p_0})\cap C([\tau,T]; H^1\cap W^{1,p }) ,  \\
\sqrt{\rho} u_t\in  L^2(0,T; L^2)\cap L^\infty(\tau,T; L^2),\quad
P_t  \in L^2(\tau,T; L^2\cap L^{p_0}),\\ \na u_t  \in L^\infty(\tau,T; L^2 )\cap    L^2(\tau,T;  L^{p_0}),\quad (\n u_t)_t \in L^2(\tau,T; L^2 ).
\end{cases}\ee
Moreover, it holds that \be \la{oiq1}\sup_{0\le t<\infty} \|\na \n\|_{L^2 }   \le 2 \|\na \n_0\|_{L^2 }  ,\quad \sup_{0\le t<\infty} \|\na \mn\|_{L^q }   \le 2 \|\na \mu(\n_0)\|_{L^q }  , \ee and that there exists some positive constant $\si$ depending only on $\|\n_0\|_{L^{3/2}}$ and $\xmu$ such that   for all $t\geq1$,
\be \la{e}
 \|\na u_t(\cdot,t)\|^2_{L^2}+
\|\na  u(\cdot,t)\|_{H^1\cap W^{1,p_0}}^2+\|P(\cdot,t)\|_{H^1\cap W^{1,p_0}}^2\le Ce^{-\sigma t},
\ee
where  $C$ depends only on $q, \beta, \bar\n, \|\n_0\|_{L^{3/2}}, \xmu, \smu, M,$ $\|\na u_0\|_{L^2},$ and  $\|\na \n_0\|_{L^2}.$
\end{theorem}

%\begin{remark}
%It should be noted here that our arguments in this paper are also valid for both the periodic boundary conditions on a torus $\mathbb{T}^3$ and the initial boundary problem in 3D bounded domain. Roughly speaking, we can reprove the global existence results in \cite{huang15,zjw} and  also derive the same large time exponential decay estimates \eqref{e} in 3D bounded domains.
%\end{remark}

As a direct consequence, our method can be applied  to the case that  $\mn\equiv\mu$ is a positive constant and obtain the following global existence and  large-time behavior of the strong solutions which improves slightly those  due to   Craig-Huang-Wang  \cite{huang13} (see   Proposition  \ref{thm3}).

\begin{theorem}\label{thm2} For constants $\bar\n>0$   and $\mu>0$,
assume that $\mn\equiv\mu$ in \eqref{1.1} and the initial data $(\rho_0, u _0)$ satisfy \eqref{hw2.2} except $u_0\in D^{2,2}.$  Then, there exists some positive constant $\ve$ depending only on $\bar \n$ such that  there exists a unique  global strong solution to the Cauchy problem \eqref{1.1} \eqref{1.3} satisfying \eqref{2.3} with $p_0=6$   provided $\|u_0\|_{\dot H^{1/2}}\le \mu\ve.$  Moreover,    it holds that
 \be \la{oiq2} \sup_{0\le t<\infty}\|\na\n\|_{L^2}\le 2\|\na\n_0\|_{L^2}, \ee and that there exists some positive constant   $\sigma$  depending only on $\|\n_0\|_{L^{3/2}}$ and $\mu$ such that  for $t\ge 1,$ \be \la{eq}
 \|\na u_t(\cdot,t)\|^2_{L^2}+
\|\na  u(\cdot,t)\|_{H^1\cap W^{1,6}}^2+\|P(\cdot,t)\|_{H^1\cap W^{1,6}}^2\le Ce^{-\sigma t} ,
\ee where   $C$ depends only on $\bar\n, \mu,$ $\|\n_0\|_{L^{3/2}}$,  $\|\na u_0   \|_{L^2},$ and $\|\na \n_0\|_{L^2}.$

\end{theorem}

A few remarks are in order.
%since the authors \cite{huang13,zhang1} show that the global strong solutions decay with algebraic rates in large time for the 3D Cauchy problem \eqref{1.1}-\eqref{1.3} .

\begin{remark}\label{re2} To the best of our knowledge,
the   exponential decay-in-time properties \eqref{e} in Theorem \ref{thm1} are
 new and somewhat surprising, since the known  corresponding decay-in-time rates for
  the strong solutions to  system \eqref{1.1}
 are algebraic   even for the constant viscosity case  \cite{huang13,zhang1}
 and  the homogeneous case \cite{chen1,sch,Kt1984,H2002}.  Moreover, as a direct consequence of \eqref{oiq1},  $\|\na \n(\cdot,t)\|_{L^2}$ remains uniformly bounded with respect to $t$   which is new even for the constant viscosity case (see \cite{huang13} or Proposition \ref{thm3}).
\end{remark}

\begin{remark} It should be noted here that our  Theorem \ref{thm1} holds
 for any function $\mn$ satisfying \eqref{n3} and for arbitrarily large initial density with vacuum (even has compact support) with  a smallness   assumption  only  on the $\dot H^\beta$-norm of the initial velocity $u_0$ with $\beta\in(1/2, 1]$, which is  in sharp contrast
to Abidi-Zhang \cite{zhang6} where they need the initial density strictly   away from vacuum and
 the smallness assumptions   on both
 $\|u_0\|_{L^2}\|\na u_0\|_{L^2}$ and $\|\mu(\n_0)-1\|_{L^\infty}$.
 %
 %
 %Furthermore, some new large time behavior of solutions are also  established in Theorem \ref{thm1}, see Remark \ref{re2} for more details.
  \end{remark}

\begin{remark}     For our case that the viscosity $\mn$ depends on $\n,$ in order to bound the $L^p$-norm of the gradient of the density, we need   the smallness conditions on the $\dot H^\beta $-norm $(\beta \in (1/2,1])$ of the initial velocity. However,  it seems that our conditions on the initial velocity may be optimal compared with the    constant viscosity case considered by  Craig-Huang-Wang \cite{huang13} where they  proved that the system \eqref{1.1} is globally wellposed for small initial data in the homogeneous Sobolev space  $\dot H^{1/2}$ which is similar to the  case of homogeneous   Navier-Stokes equations (see \cite{ft1964}).   Note that for the case of   initial-boundary-value problem in 3D bounded domains, Huang-Wang \cite{huang15} and Zhang  \cite{zjw}  impose smallness conditions on  $\|\na u_0\|_{L^2}.$ Furthermore, in our Theorems 1.2 and 1.3, there is no need to imposed  additional initial compatibility conditions, which is assumed in \cite{huang13,huang15,zjw} for the global existence of strong solutions.
\end{remark}

\begin{remark}It is easy to prove that  the strong-weak uniqueness theorem   \cite[Theorem 2.7]{L1996}  still holds for the initial data $(\n_0,u_0)$ satisfying \eqref{2.2}  after modifying its proof slightly. Therefore, our
Theorem  \ref{thm1}  can be regarded as the uniqueness and
 regularity
theory of Lions's weak solutions \cite{L1996} with the initial velocity suitably small in the
 $ {\dot H^\beta}$-norm.
\end{remark}

\begin{remark}    In \cite{kim2004}, Cho-Kim considered the initial boundary value problem in 3D bounded smooth domains.  In addition  to \eqref{2.2}, assuming   that the initial data satisfy  the following  compatibility conditions  \bnn \la{zxq1}
- \div(\mu(\n_0)(\na u_0+(\na u_0)^{\rm T}))  +\na P_0=\n_0^{1/2}g
\enn for some $(P_0,g)\in H^1\times L^2,$ it is shown (\cite{kim2004}) that  the local-in-time strong solution $(\n,u)$ satisfies \be \la{cont1}\n u_t\in C([0,T];L^2). \ee  However,  to  obtain \eqref{cont1}, it seems difficult to follow the proof of \eqref{cont1} as in  \cite{kim2004}. Indeed, in our Proposition \ref{lem5.a3},
%without the condition  \eqref{zxq1},
we give a complete new proof to show that  $\n u_t\in H^1(\tau,T;L^2)$ (for any $0<\tau<T<\infty$) which directly implies \be \la{cont2}\n u_t\in C([\tau,T];L^2). \ee In fact, \eqref{cont2} is crucial for deriving  the time-continuity of $\na u$ and $P,$  that is (see \eqref{2.3}),  \be \na u, P\in C([\tau,T];H^1\cap W^{1,p}).\ee
\end{remark}

We now make some comments on  the analysis in this paper. To extend the local  strong solutions whose existence is obtained by Lemma \ref{local} globally in time,
  one needs  to establish global a priori estimates on smooth solutions   to \eqref{1.1}--\eqref{1.3}  in suitable higher norms. It turns out  that as in the 3D bounded case \cite{huang15,zjw}, the   key ingredient  here is to get the time-independent bounds on the $L^1(0,T; L^\infty)$-norm of $\na u$ and then the $L^\infty(0,T; L^q)$-norm of $\na \mn$  and the $L^\infty(0,T; L^2)$-one of $\na \n$.
However, as mentioned by Huang-Wang \cite{huang15},
the methods used in \cite{huang15,zjw}  depend  crucially on the boundedness of the domains.  Hence, some new ideas  are needed here. First, using the initial layer
analysis (see  \cite{Hoh4,HLX2012}) and an interpolation argument (see  \cite{stein}),  we   succeed  in bounding the $L^1(0,\min\{1,T\}; L^\infty)$-norm of $\na u$    by $\|u_0\|_{\dot H^\beta}$ (see \eqref{3d4.8}).  Then,  in order  to estimate the $L^1(\min\{1,T\},T; L^\infty)$-norm of $\na u$, we find that  $\|\n^{1/2}u(\cdot,t)\|^2_{L^2}$ in fact decays at the rate of $e^{-\si t} (\si >0)$ for large time (see \eqref{ed1}), which can be achieved by combining the standard energy equality  (see \eqref{ed1.1}) with the following fact
$$
 \|\n^{1/2}u\|_{L^2}^2\le \|\n\|_{L^{3/2}}\|u\|_{L^6}^2\le C\|\na u\|_{L^2}^2,$$ due to \eqref{1.1}$_1$ and the Sobolev inequality. With this key exponential decay-in-time rate at hand, we can obtain that both $\|\na u(\cdot,t)\|^2_{L^2}$ and  $\|\n^{1/2}u_t(\cdot,t)\|^2_{L^2}$ decay  at the same rate as $e^{-\si t} (\si >0)$ for large time (see \eqref{ed2} and \eqref{ed3}). In fact,   all these  exponential decay-in-time rates are the key to   obtaining  the desired  uniform  bound (with respect to time) on the  $L^1(\min\{1,T\},T; L^\infty)$-norm
of $\na u$   (see \eqref{jia9}). Finally, using these a priori estimates and the fact that the velocity is divergent free, we   establish the time-independent estimates  on
the gradients of the density and the velocity which guarantee  the extension of local strong solutions (see Proposition \ref{lem5.a3}).

The rest of this paper is organized as follows. In Section \ref{sec2}, we collect some elementary facts and inequalities that will be used later. Section \ref{sec3} is devoted to the a priori estimates. Finally, we will prove  Theorems  \ref{thm1} and \ref{thm2} in Section 4.

\section{Preliminaries}\label{sec2}
In this section we shall enumerate some auxiliary lemmas.

We start with the local existence of strong solutions which has  been proved   in \cite{hls}.
\begin{lemma}\label{local}
Assume that $(\rho_0, u_0)$ satisfies   \eqref{2.2} except $u_0\in\dot H^\beta.$  Then there exist  a small time $T_0>0$ and a unique strong solution $(\rho, u, P)$ to the problem \eqref{1.1}--\eqref{1.3} in $\mathbb{R}^{3}\times(0,T_0)$ satisfying \eqref{2.3}.
\end{lemma}

% Considering the following Stokes type equations
% \be\ba\label{3rd1}
%\begin{cases}
% -\div(\mn \na u)  +\nabla P=F, \\
%\div u=0.
%\end{cases}
%\ea\ee

The following regularity results on the Stokes equations will be useful for our derivation of higher order a priori estimates.
\begin{lemma} \la{stokes} For   positive constants $\xmu,\smu,$ and $ q\in (3 ,\infty)$, in addition to \eqref{n3},  assume that    $\mn$ satisfies \be\la{ij80}\na\mn\in L^q,\quad
  0<\xmu\le \mn\le\smu<\infty. \ee
   Then, if $F\in L^{6/5}\cap L^r$ with   $r\in[ 2q/(q+2),q],$ there exists some positive constant $C$ depending only on $ \xmu , \smu , r, $ and $q$  such that the unique weak solution $(u,P)\in D^1_{0,\sigma}\times L^2$ to the following Cauchy problem
\be\label{3rd1}
\begin{cases}
 -\div(2\mn d)  +\nabla P=F,\,\,\,\,&x\in \rrr,\\
 \div u=0,   \,\,\,&x\in  \rrr,\\
u(x)\rightarrow0,\,\,\,\,&|x|\rightarrow\infty,
\end{cases}
\ee satisfies
 \be\ba\label{3rd2}
\|\na u\|_{L^2 }+\|P\|_{L^2 }\le C \|F\|_{L^{6/5} },
\ea\ee
 \be\ba\label{3rd3'}
\|\na^2 u\|_{L^r}+\xl\|\na P\xr\|_{L^r}\le  C  \|F\|_{L^r}+ C \|\na\mn\|_{L^q}^{\frac{q(5r-6)}{2r(q-3)}}\|F\|_{L^{6/5}} .
\ea\ee
Moreover, if $F=\div g$ with $g\in L^2\cap L^{\ti r}$ for some $\ti r\in (6q/(q+6),q],$ there exists a positive constant $C$
 depending only on $\xmu, \smu, q,$ and $\ti r$ such that the unique weak solution $(u,P)\in D^1_{0,\sigma}\times L^2$ to
 \eqref{3rd1}  satisfies\be \la{3e1}\|\na u\|_{L^2\cap L^{\ti r}}+ \|P\|_{L^2\cap L^{\ti r}}\le  C  \|g\|_{L^2\cap L^{\ti r} }+C \|\na\mn\|_{L^q}^{\frac{3q(\ti r-2)}{2\ti r(q-3)}}\|g\|_{L^2}.\ee\end{lemma}

\pf First, multiplying \eqref{3rd1}$_1$ by $u$ and integrating by parts, we obtain after using \eqref{3rd1}$_2$ that
  \be\ba\label{3rd3}\notag
 2\int \mn|d|^2dx=\int F\cdot udx\le \|F\|_{L^{6/5}}\|u\|_{L^6}\le C\|F\|_{L^{6/5}} \|\na  u\|_{L^2},
\ea\ee
which together with \eqref{ij80} yields
  \be\ba\label{3rd4}
\|\na  u\|_{L^2} \le C\xmu^{-1} \|F\|_{L^{6/5}},
\ea\ee due to
  \be\ba\label{3rd12}
 2\int|d|^2dx= \int|\na u|^2dx.
\ea\ee

Furthermore, it follows from \eqref{3rd1}$_1$ that
 \be\ba\label{3rd5}\notag
 P=-(-\Delta)^{-1}\div F-(-\Delta)^{-1}\div\div(2\mn d),
\ea\ee
which together with the Sobolev inequality and \eqref{3rd12}   gives
 \be\ba\label{3rd6}\notag
 \|P\|_{L^2} &\le \|(-\Delta)^{-1}\div F\|_{L^2} +\|2\mn d\|_{L^2} \le C\|F\|_{L^{6/5}} +C \| \na u\|_{L^2} .
\ea\ee
  Combining this with \eqref{3rd4} leads to \eqref{3rd2}.

Next, we rewrite  \eqref{3rd1}$_1$  as
 \be\ba\label{3rd7}
  -\Delta u+\na \xl(\frac{P}{\mn}\xr)=\frac{F}{\mn}+\frac{2d\cdot\na\mn }{\mn}-\frac{P\na\mn}{\mn^2}.
\ea\ee
Applying the standard $L^p$-estimates to the Stokes system \eqref{3rd7} \eqref{3rd1}$_2$ \eqref{3rd1}$_3$  yields that for  $r\in[ 2q/(q+2),q],$
\bnn\ba\label{3rd8}
  \|\na^2 u\|_{L^r}+\xl\|\na P\xr\|_{L^r}
 &\le \|\na^2 u\|_{L^r}+C\xl\|\na  \xl(\frac{P}{\mn}\xr)\xr\|_{L^r}+C\xl\| \frac{P\na\mn}{\mn^2}\xr\|_{L^r}\\
 &\le C \xl\|  {F}  \xr\|_{L^r}+C\xl\|  {2d\cdot\na\mn} \xr\|_{L^r}+C\xl\|  {P\na\mn} \xr\|_{L^r}  \\
   &\le C \|F\|_{L^r}+C \|\na\mn\|_{L^q} \|\na u\|_{L^2}^{\frac{2r(q-3)}{q(5r-6)}}\|\na^2 u\|_{L^r}^{1-\frac{2r(q-3)}{q(5r-6)}}\\
   &\quad+C \|\na\mn\|_{L^q}\xl\| P\xr\|_{L^2}^{\frac{2r(q-3)}{q(5r-6)}}\xl\|\na P \xr\|_{L^r}^{1-\frac{2r(q-3)}{q(5r-6)}}\\
 % &\le C \xmu^{-1}\|F\|_{L^r}+C\xmu^{-1}\|\na\mn\|_{L^q}\|\na u\|_{L^2}^{a_r}\|\na^2 u\|_{L^r}^{1-a_r}\\
%  &\quad+C\xmu^{-1}\|\na\mn\|_{L^q}\xl\| \frac{P}{\mn}\xr\|_{L^2}^{a_r}\xl\|\na \xl(\frac{P}{\mn}\xr)\xr\|_{L^r}^{1-a_r} \\
    &\le C \|F\|_{L^r}+C \|\na\mn\|_{L^q}^{\frac{q(5r-6)}{2r(q-3)}} (\|\na u\|_{L^2}+  \|  {P} \|_{L^2} ) \\&\quad+\frac{1}{2}\xl(\|\na^2 u\|_{L^r}+ \xl\|\na P\xr\|_{L^r}\xr),
\ea\enn
  which combined with  \eqref{3rd2} yields    \eqref{3rd3'}.
% \be\ba\label{3rd9}
% &\|\na^2 u\|_{L^r}+\xl\|\na \xl(\frac{P}{\mn}\xr)\xr\|_{L^r}\\
% &\le C \xmu^{-1}\|F\|_{L^r}+C\xmu^{-1-a_r^{-1}}\xl(1+\smu\xmu^{-1}\xr)\|\na\mn\|_{L^q}^{a_r^{-1}}\|F\|_{L^{6/5}}.
%\ea\ee

Finally, we will prove \eqref{3e1}. Multiplying \eqref{3rd1}$_1$ by $u$ and integrating by parts lead to \bnn 4\int\mu(\n)|d|^2dx=-2\int g\cdot\na udx\le  {\xmu} \|\na u\|_{L^2}^2+C\|g\|_{L^2}^2,\enn which together with \eqref{3rd12}  gives \be \la{5rd1}\|\na u\|_{L^2}\le C\|g\|_{L^2}.\ee It follows from \eqref{3rd1}$_1$ that
 \bnn\ba
 P=-(-\Delta)^{-1}\div \div g-(-\Delta)^{-1}\div\div(2\mn d),
\ea\enn which  implies that for any $p\in [2,\ti r],$\be \la{5r2t}\|P\|_{L^p}\le C(p)\|\na u\|_{L^p}+ C(p)\|g\|_{L^p}.\ee
In particular,   this   combined with \eqref{5rd1} shows \be \la{5r2}\|P\|_{L^2}+\|\na u\|_{L^2}\le C\|g\|_{L^2}.\ee

Next, we rewrite  \eqref{3rd1}$_1$  as
 \be\ba\label{5rd7}
  -\Delta u+\na \xl(\frac{P}{\mn}\xr)=
  \div\left(\frac{g}{\mn}\right) +\ti G,
\ea\ee where $$\ti G\triangleq\frac{g\cdot\na\mn }{\mn^2}+\frac{2d\cdot\na\mn }{\mn}-\frac{P\na\mn}{\mn^2} $$
satisfies for any $\ve>0,$ \be\la{5r3}\ba \|\ti G\|_{L^{\frac{3\ti r}{ 3+\ti r}}}&\le \ve (\|g\|_{L^{\ti r}}+\|\na u\|_{L^{\ti r}}+\|P\|_{L^{\ti r}})\\&\quad+C(\ve)\|\na\mn\|_{L^q}^{\frac{3q(\ti r-2)}{2\ti r(q-3)}}( \|g\|_{L^2}+\|\na u\|_{L^2}+\|P\|_{L^2}).\ea\ee
Using \eqref{5rd7} and \eqref{3rd1}$_3$, we have\bnn\ba \|\na u\|_{L^{\ti r}}&\le C\|\na\times u\|_{L^{\ti r}}\\ &=C\|(-\Delta)^{-1}\na \times \div (g(\mn)^{-1}) +(-\Delta)^{-1}\na \times \ti G\|_{L^{\ti r}}\\ &\le C\|g\|_{L^{\ti r}}+C\|\ti G\|_{L^{\frac{3\ti r}{3+\ti r}}},\ea\enn which together with \eqref{5r2t} yields  \bnn\ba \|\na u\|_{L^{\ti r}}+\|P\|_{L^{\ti r}} \le C\|g\|_{L^{\ti r}}+C\|\ti G\|_{L^{\frac{3\ti r}{3+\ti r}}}.\ea\enn  Combining  this,   \eqref{5r3}, and  \eqref{5r2} gives \eqref{3e1}. The proof of Lemma \ref{stokes} is  finished. \hfill $\Box$

\section{A Priori Estimates}\label{sec3}
In this section, we will establish some necessary a priori bounds of local strong solutions $(\rho,u,P)$ to the Cauchy problem \eqref{1.1}--\eqref{1.3}  whose existence is guaranteed by Lemma \ref{local}. Thus, let $T>0$ be a fixed time and $(\rho, u,P)$  be the smooth solution to \eqref{1.1}-\eqref{1.3} on $\rrr\times(0,T]$ with smooth initial data $(\rho_0,u_0)$ satisfying \eqref{2.2}.

We have the following key a priori estimates on $(\n,u,P)$.
\begin{proposition}\la{pr1}  There exists some  positive constant  $\ve_0$
    depending    only on  $q, \beta, \bar\n, \xmu, \smu, $ $\|\n_0\|_{L^{3/2}}, $ and $M$  such that if
       $(\n,u,P)$  is a smooth solution of
       \eqref{1.1}--\eqref{1.3}  on $\rrr \times (0,T] $
        satisfying
 \be\la{3a1}
     \sup_{t\in[0,T]}\|\na\mn\|_{L^q}\le 4M ,\quad  \int_0^T\|\na u\|_{L^2}^4dt \leq 2\|u_0\|_{\dot H^\beta}^2 ,
   \ee
 the following estimates hold
        \be\la{3a2}
\sup_{t\in[0,T]}\|\na\mn\|_{L^q} \le 2M ,\quad \int_0^T\|\na u\|_{L^2}^4dt \leq  \|u_0\|_{\dot H^\beta}^2 ,
  \ee
   provided $\|u_0\|_{\dot H^\beta}\le \ve_0.$
\end{proposition}

Before proving Proposition \ref{pr1}, we establish some necessary a priori estimates, see Lemmas \ref{lem3.3}--\ref{lem5.1}.

%Moreover, we assume that the following a priori hypothesis holds for the viscosity $\mn$:
%\be\label{3a1}\sup_{t\in[0,T]}\|\na\mn\|_{L^q}\le 4M,~~\mbox{for}~q>3,~~M\triangleq\|\na \mu(\n_0)\|_{L^q}.\ee

We start with the following   time-weighted  estimates on the $L^\infty(0,\min\{1,T\};L^2)$-norm of the gradient of  velocity.
\begin{lemma}\label{lem3.3} Let $(\n, u, P)$ be a  smooth  solution to  \eqref{1.1}--\eqref{1.3} satisfying \eqref{3a1}. Then   there exists  a generic positive constant  $C$  depending only on  $q,$ $\beta,$ $\ban,$ $\xmu,$ $\smu,$ $\|\n_0\|_{L^{3/2}},$  and $M$ such that
\be\label{gj1} \sup_{t\in[0,{\zeta(T)}]}\xl(t^{1-\beta}\|\na u\|_{L^2}^2\xr)+\int_0^{\zeta(T)} t^{1-\beta}\|\n^{1/2} u_t\|_{L^2}^2dt\le C \|u_0\|_{\dot H^\beta}^2,\ee where $\zeta(t)$ is defined by  \bnn {\zeta(t)} \triangleq \min\{1,t\}.\enn
\end{lemma}

{\it Proof.} First, standard arguments (\cite{L1996}) imply that \begin{equation}\label{gj0}
0\le\n\leq \bar\n,~~~~~\|\n\|_{L^{3/2}}= \|\n_0\|_{L^{3/2}}.
\end{equation}

%Next

Next, for fixed $(\n, u)$, we  consider  the following linear  Cauchy problem for $(w,\tilde P)$:
\be\label{z1.1}
\begin{cases}
 \n w_t+\n u\cdot\na w-\div\xl(\mn\xl[\na w+(\na w)^T\xr] \xr)  +\nabla \tilde P= 0,\, &x\in \rrr,\\
 \div w=0,   \, &x\in  \rrr,\\
 w(x,0)=w_0,   \,  &x\in  \rrr.
%w(x,\cdot)\rightarrow0,\,\,\,\,&|x|\rightarrow\infty,
\end{cases}
\ee
It follows from  Lemma \ref{stokes}, \eqref{z1.1}$_1$, \eqref{3a1},   \eqref{gj0}, and the Garliardo-Nirenberg inequality that
\bnn%\label{3d4}
\ba
 \|\na  w\|_{H^1}+ \|  \tilde P \|_{H^1}
 %&\le C \frac{1}{\xmu} \xl(\|\tilde F\|_{L^2}+ \frac{\xmu+\smu}{\xmu^{1+a_2^{-1}}}\|\na\mn\|_{L^q}^{a_2^{-1}}\|\tilde F\|_{L^{6/5}}\xr)\\
% &\le C\xl(\ban^{1/2}\|\n^{1/2} w_t+\n^{1/2} u\cdot\na w\|_{L^2}+\|\n^{1/2}\|_{L^3}\|\n^{1/2} w_t+\n^{1/2} u\cdot\na w\|_{L^2}\xr)\\
  &\le C \xl(\|\n  w_t+\n u\cdot\na w\|_{L^2}+\|\n  w_t+\n u\cdot\na w \|_{L^{6/5}}\xr)\\
   &\le C ( \ban^{1/2}+\|\n\|^{1/2}_{L^{3/2}})\xl(\|\n^{1/2} w_t\|_{L^2}+\ban^{1/2} \|  u\cdot \na  w\|_{L^2}\xr) \\
  %&\le C(q,\ban,\tn,\xmu,\smu,M) \xl(\|\n^{1/2} u_t\|_{L^2}+\|\n^{1/2}u\cdot\na u\|_{L^2}\xr)\\
   %&\le C(q,\ban,\tn,\xmu,\smu,M) \|\n^{1/2} u_t\|_{L^2}+C(q,\ban,\tn,\xmu,\smu,M)\|u\|_{L^6}\|\na u\|_{L^3}\\
      &\le C  \|\n^{1/2} w_t\|_{L^2}+C \|\na u\|_{L^2} \|\na w\|_{L^2}^{1/2}\|\na^2 w\|_{L^2}^{1/2}\\
      &\le C \|\n^{1/2} w_t\|_{L^2}+C \|\na u\|_{L^2}^2 \|\na w\|_{L^2}+\frac12 \|\na^2 w\|_{L^2},
\ea\enn which directly yields that
\begin{equation} \label{3d4}
\ba
 &\|\na w\|_{H^1}+ \|  \tilde P \|_{H^1}+\|\n   w_t+\n u\cdot\na w \|_{L^{6/5}\cap L^2}
 %&\le C \frac{1}{\xmu} \xl(\|\tilde F\|_{L^2}+ \frac{\xmu+\smu}{\xmu^{1+a_2^{-1}}}\|\na\mn\|_{L^q}^{a_2^{-1}}\|\tilde F\|_{L^{6/5}}\xr)\\
% &\le C\xl(\ban^{1/2}\|\n^{1/2} w_t+\n^{1/2} u\cdot\na w\|_{L^2}+\|\n^{1/2}\|_{L^3}\|\n^{1/2} w_t+\n^{1/2} u\cdot\na w\|_{L^2}\xr)\\
  \\      &\le C \|\n^{1/2} w_t\|_{L^2}+C \|\na u\|_{L^2}^2 \|\na w\|_{L^2} .
\ea\end{equation}
Multiplying \eqref{z1.1}$_1$ by $w_t$ and  integrating the resulting equality  by parts lead to
\begin{equation}\label{z1.3}%\label{3d3}
\ba
& \frac{1}{4}\frac{d}{dt}\int \mn \xl| \na w+(\na w)^T\xr|^2dx+\int\rho|w_t|^{2}dx\\
&=-\int\n u\cdot\na w\cdot w_tdx+ \frac{1}{4}\int  \mn u\cdot\na \xl|\na w+(\na w)^T\xr|^2dx\\
&\le \ban^{1/2}\|\n^{1/2}w_t\|_{L^2}\|u\|_{L^6}\|\na w\|_{L^3}+C\smu \|u\|_{L^6}\|\na w\|_{L^3}\|\na^2 w\|_{L^2} \\
&\le C \|\n^{1/2}w_t\|_{L^2}  \|\na u\|_{L^2}  \|\na w\|_{L^2}^{1/2} \|\na^2 w\|_{L^2}^{1/2}   +C  \|\na u\|_{L^2} \|\na w\|_{L^2}^{1/2}\|\na^2 w\|_{L^2}^{3/2} \\
&\le \frac{3}{4}\|\n^{1/2}w_t\|_{L^2}^2+C \|\na u\|_{L^2}^4 \|\na w\|_{L^2}^2,
\ea\end{equation}
 where in the last inequality one has used   \eqref{3d4}. This  combined with Gr\"onwall's inequality  and \eqref{3a1}  yields
\begin{equation}\label{3d6}\ba
 \sup_{t\in[0,{\zeta(T)}]}\int |\na w|^2dx+\int_0^{\zeta(T)}\int\rho|w_t|^{2}dxdt
 \le C \|\na w_0\|_{L^2}^2.
\ea\end{equation}

Furthermore, multiplying \eqref{z1.3} by $t$ leads to
\begin{equation}\notag%\label{3d7}
\ba
& \frac{d}{dt}\xl(t\int  \mn \xl| \na w+(\na w)^T\xr|^2dx\xr)+ t\int\rho|w_t|^{2}dx\\&
 \le C  t\|\na w\|_{L^2}^2\|\na u\|_{L^2}^4+C\|\na w\|_{L^2}^2.
\ea\end{equation}
Combining this  with Gr\"onwall's inequality and \eqref{3a1} shows
\begin{equation}\label{3d8}\ba
 \sup_{t\in[0,{\zeta(T)}]}t\int |\na w|^2dx+\int_0^{\zeta(T)} t\int\rho|w_t|^{2}dxdt
%&\le C  \int \rho_0  |w_0|^2dx\\
&\le C  \|w_0\|_{L^2}^2,
\ea\end{equation} where one has used the following simple fact
\bnn\ba\label{z1.2}
 \sup_{t\in[0,{\zeta(T)}]}\|\rho^{1/2} w \|_{L^2}^2
+ \int_{0}^{{\zeta(T)}}\|\nabla w \|_{L^2}^{2}dt
& \le  C\|w_0\|_{L^2}^2,
\ea\enn which can be obtained by
multiplying \eqref{z1.1}$_1$ by $w$ and integrating by parts.

  Hence,
the standard Stein-Weiss interpolation arguments (see \cite{stein})
 together with \eqref{3d6} and \eqref{3d8} imply that for any $\theta\in[\beta,1]$,
\begin{equation} \label{3d9}
\ba
&\sup_{t\in[0,{\zeta(T)}]}t^{1-\theta}\int |\na w|^2dx+\int_0^{\zeta(T)} t^{1-\theta}\int\rho|w_t|^{2}dxdt \le C(\theta) \|w_0\|_{\dot H^\theta}^2.
\ea\end{equation}
%for some positive constant $C_2=C_2(q,,\beta\ban,\tn,\xmu,\smu,M)$.

Finally, taking $w_0=u_0$, the uniqueness of strong solutions to the linear problem \eqref{z1.1}  implies that $w\equiv u.$   The estimate \eqref{gj1} thus follows from \eqref{3d9}. The proof of Lemma \ref{lem3.3} is finished.  \hfill $\Box$

As an application of Lemma \ref{lem3.3}, we   have   the following time-weighted  estimates on  $\|\n^{1/2} u_t\|_{L^2}^2 $ for small time.

\begin{lemma}\label{lem3.4}Let $(\n, u, P)$ be a  smooth  solution to  \eqref{1.1}--\eqref{1.3} satisfying \eqref{3a1}. Then  there exists  a generic positive constant  $C$  depending only on  $q,$ $\beta,$ $\ban,$ $\xmu,$ $\smu,$ $\|\n_0\|_{L^{3/2}},$ and $M$ such that
\be\label{gj2} \sup_{t\in[0,{\zeta(T)}]}\xl(t^{2-\beta}\|\n^{1/2} u_t\|_{L^2}^2\xr)+\int_0^{\zeta(T)} t^{2-\beta}\|\na  u_t\|_{L^2}^2dt\le C \|u_0\|_{\dot H^\beta}^2.\ee
\end{lemma}

{\it Proof.}  First,
operating $\partial_{t}$ to \eqref{1.1}$_2$ yields that
\be \la{utt}\ba &
\n u_{tt}+ \n u \cdot\na u_t-\div(2\mn d_t)+\na P_t\\&=-\n_tu_t-(\n u)_t\cdot\na u+\div(2(\mn)_t d).\ea
\ee
Multiplying the above equality by $u_t$, we obtain after using integration by parts and \eqref{1.1}$_1$ that
\be\ba\label{3d2.1}&
\frac12\frac{d}{dt}\int\rho|u_t|^{2}dx+\int2\mn|d_t|^{2}dx   \\
%&=-\int \n_t|u_t|^2dx-\int\n u\cdot\na(u\cdot\na u\cdot u_t)dx\\
%&\quad-\int\n u_t\cdot\na u\cdot u_tdx-2\int (\mn)_t\cdot d\cdot\na u_tdx\\
&=-2\int \n u\cdot\na u_t\cdot u_tdx-\int\n u\cdot\na(u\cdot\na u\cdot u_t)dx\\
&\quad-\int\n u_t\cdot\na u\cdot u_tdx+2\int \xl(u\cdot\na\mn\xr) d\cdot\na u_tdx \triangleq \sum_{i=1}^4J_i.
\ea\ee

Now, we will use the Gagliardo-Nirenberg inequality, \eqref{3a1}, and \eqref{gj0} to estimate each term  on the right hand of \eqref{3d2.1} as follows:
\be\ba\label{3d2.2}
 | J_1|+|J_3|
  &\le C \|\n^{1/2}u_t\|_{L^3}\|\na u_t\|_{L^2}\|u\|_{L^6}+C \|\n^{1/2}u_t\|_{L^3}\|\na u\|_{L^2}\|u_t\|_{L^6} \\
  &\le C \|\n^{1/2}u_t\|_{L^2}^{1/2}\|\na u_t\|_{L^2}^{3/2}\|\na u \|_{L^2}\\
    &\le  \frac{1}{4}\xmu\|\na u_t\|_{L^2}^2+C \|\n^{1/2}u_t\|_{L^2}^2\|\na u \|_{L^2}^4,
\ea\ee
\be\ba\label{3d2.3}
  |J_2| &=\left|\int\n u\cdot\na(u\cdot\na u\cdot u_t)dx\right|\\
  &\le C\int \n|u||u_t|\left(|\na u|^2+|u| |\na^2u|\right)dx+\int\n|u|^2|\na u||\na u_t|dx\\
  &\le C\|u\|_{L^6}\|u_t\|_{L^6} \xl( \|\na u\|_{L^3}^2+   \|u\|_{L^6} \|\na^2 u\|_{L^2} \xr)+ C  \|u\|_{L^6}^2\|\na u\|_{L^6}\|\na u_t\|_{L^2}\\
  &\le C  \|\na u_t\|_{L^2}\|\na^2 u\|_{L^2}\|\na u\|_{L^2}^{2}\\
    &\le  \frac{1}{8}\xmu\|\na u_t\|_{L^2}^2+C \|\na^2 u\|_{L^2}^2 \|\na u\|_{L^2}^4,
\ea\ee and
\be\ba\label{3d2.5}
  |J_4 |   &\le  C\|\na\mn\|_{L^q} \|u\|_{L^\infty} \|\na u_t\|_{L^2}\|\na u\|_{L^{\frac{2q}{q-2}}}\\    &\le  C(q,M)\|u\|_{L^6}^{1/2}\|\na u\|_{L^6}^{1/2} \|\na u_t\|_{L^2}\|\na u\|_{L^2}^{\frac{q-3}{q}}\|\na^2 u\|_{L^2}^\frac{3}{q} \\    &\le  \frac{1}{8}\xmu\|\na u_t\|_{L^2}^2+C \|\na u \|_{L^2} \|\na^2 u\|_{L^2}^3+ C \|\na u \|_{L^2}^4.
\ea\ee

  Substituting \eqref{3d2.2}--\eqref{3d2.5} into \eqref{3d2.1}  gives
\be\ba\label{3d2.6}&
 \frac{d}{dt}\int\rho|u_t|^{2}dx+\xmu\int|\nabla u_t|^{2}dx   \\
&\le  C \left(\|\rho^{1/2} u_t\|_{L^2}^2 +   \|\na^2 u\|_{L^2}^2 \right)\|\na u\|_{L^2}^4+C \|\na u \|_{L^2} \|\na^2 u\|_{L^2}^3+ C \|\na u \|_{L^2}^4 \\
&\le C \|\rho^{1/2}u_t\|_{L^2}^2 \|\na u \|_{L^2}^4 +C \|\rho^{1/2}u_t\|_{L^2}^3\|\na u \|_{L^2}  +C \|\na u \|_{L^2}^{10} +C \|\na u \|_{L^2}^{2},
\ea\ee
where in the last inequality one has used
  \be\label{z1.5}
\ba
 &\|\na  u\|_{H^1}+ \| P\|_{H^1}+ \|\n (u_t+u\cdot \na u)\|_{L^{6/5}\cap L^2}  \\& \le C  \xl(\|\n^{1/2} u_t\|_{L^2}+ \|\na u\|_{L^2}^3 \xr),
\ea\end{equation} which can be obtained by taking $w\equiv u$ in \eqref{3d4}.
  It thus follows   from \eqref{3d2.6}  and \eqref{gj1} that for $ t \in (0,{\zeta(T)}]$, \be\ba\label{3d2.6'}&
 \frac{d}{dt}\int\rho|u_t|^{2}dx+\xmu\int|\nabla u_t|^{2}dx
\\
&\le C \|\rho^{1/2}u_t\|_{L^2}^2 \xl(\|\na u \|_{L^2}^4 + \|\rho^{1/2}u_t\|_{L^2} \|\na u \|_{L^2}\xr)\\&\quad +C  t^{3(\beta-1)} \|\na u \|_{L^2}^{4}+C \|\na u \|_{L^2}^{2}.
\ea\ee Since  \eqref{gj1} implies
\bnn\label{zcd2}\ba
 & \int_0^{\zeta(T)} \|\rho^{1/2} u_t\|_{L^2}\|\na u\|_{L^2}dt\\&\le C\sup_{0\le t\le {\zeta(T)}}  \xl( t^{ \frac{1-\beta}{2}}\|\na u\|_{L^2}\xr) \xl( \int_0^{\zeta(T)} t^{1-\beta}\|\sqrt{\rho} u_t\|_{L^2} ^2 dt\xr)^{1/2} \xl( \int_0^{\zeta(T)} t^{ 2 \beta-2}dt\xr)^{1/2} \\&  \le C \|u_0\|_{\dot H^\beta}^2 ,
\ea\enn
  we multiply \eqref{3d2.6'} by $t^{2-\beta} $   and use Gr\"onwall's inequality,   \eqref{3a1}, and \eqref{gj1}  to obtain \eqref{gj2}.  The proof of Lemma \ref{lem3.4}  is  finished. \hfill $\Box$

Next,   we will prove the following      exponential  decay-in-time  estimates on the solutions for large time, which   plays a crucial role  in our analysis.

\begin{lemma}\label{lem-ed1}Let $(\n, u, P)$ be a smooth solution to  \eqref{1.1}--\eqref{1.3} satisfying \eqref{3a1}.
Then for   \be \la{laj1}\sigma\triangleq    3\xmu/(4\|\n_0\|_{L^{3/2}}) , \ee  there exists  a generic positive constant  $C$  depending only on  $q,$ $\beta,$ $\ban,$  $\xmu,$ $\smu,$ $\|\n_0\|_{L^{3/2}},$  and $M$ such that
%it holds
\be\label{ed1} \sup_{t\in[0,T]}e^{\sigma t}\|\n^{1/2}u\|_{L^2}^2 +\int_0^Te^{\sigma t}\int |\na u|^2dxdt \le  C \|u_0\|_{\dot H^{\beta}}^2,\ee
\be\label{ed2}   \sup_{t\in[{\zeta(T)},T]}e^{\sigma t}\int |\na u|^2dx +\int_{\zeta(T)}^T e^{\sigma t}\int\rho|u_t|^{2}dxdt \le C  \|u_0\|_{\dot H^{\beta}}^2, \ee
\be\label{ed3}   \sup_{t\in[{\zeta(T)},T]}e^{\sigma t}\int\rho|u_t|^{2}dx +\int_{\zeta(T)}^T e^{\sigma t}\int  |\na u_t|^2dx dt \le C  \|u_0\|_{\dot H^{\beta}}^2, \ee
 and
\begin{equation}\label{ed4}\ba
  \sup_{t\in[{\zeta(T)},T]}e^{\sigma t}\xl(\|\na  u\|_{H^1}^2+\|P\|_{H^1}^2\xr) \le  C \|u_0\|_{\dot H^\beta}^2.
\ea\end{equation}
\end{lemma}

{\it Proof.}
First, multiplying \eqref{1.1}$_2$ by $u$ and integrating by parts lead to
\be\ba\label{ed1.1}
 \frac{1}{2}\frac{d}{dt}\|\n^{1/2}u\|_{L^2}^2+\int 2\mn |d|^2dx=0.
\ea\ee
It follows from the Sobolev inequality  \cite[(II.3.11)]{N1959},   \eqref{gj0}, and \eqref{3rd12} that
\be\ba\label{ed1.2}
 \|\n^{1/2}u\|_{L^2}^2\le \|\n\|_{L^{3/2}}\|u\|_{L^6}^2\le  \frac{4}{3}\|\n_0\|_{L^{3/2}}\|\na u\|_{L^2}^2\le    \si^{-1}\int 2\mn |d|^2dx,
\ea\ee
with  $\si$ is defined as in \eqref{laj1}. Putting \eqref{ed1.2} into  \eqref{ed1.1}     yields
$$  \frac{d}{dt}\|\n^{1/2}u\|_{L^2}^2+\si \|\n^{1/2}u\|_{L^2}^2+\int 2\mn |d|^2dx\le 0,
$$ which together with Gr\"onwall's inequality gives
\be\ba\label{3d2}&
 \sup_{t\in[0,T]}e^{\sigma t}\|\n^{1/2}u\|_{L^2}^2 +\int_0^Te^{\sigma t}\int |\na u|^2dxdt\\ &\le C \|\n_0^{1/2}u_0\|_{L^2}^2 \le C\|\n_0\|_{L^{\frac{3}{2\beta}}}\|u_0\|_{L^{\frac{6}{3-2\beta}}}^2  \le  C \|u_0\|_{\dot H^{\beta}}^2,
\ea\ee
  due to $\beta\in(1/2,1]$.

Next, similar to \eqref{z1.3}, we have
\be\label{3d5}
\ba
& \frac{d}{dt}\int  2\mn  |d|^2dx+\int\rho|u_t|^{2}dx\le C  \|\na u\|_{L^2}^4\|\na u\|_{L^2}^2,
\ea\ee
which   combined with Gr\"onwall's inequality, \eqref{3d2}, \eqref{gj1},  and \eqref{3a1} gives \eqref{ed2}.
%\begin{equation}\label{ed2.2}\ba
%  \sup_{t\in[1,T]}e^{\sigma t}\int   |\na u|^2dx +\int_1^T e^{\sigma t}\int\rho|u_t|^{2}dxdt \le C  \|u_0\|_{\dot H^{\beta}}^2.
%\ea\end{equation}

Furthermore, multiplying \eqref{3d2.6} by $e^{\sigma t},$ we obtain \eqref{ed3} after using Gr\"onwall's inequality, \eqref{gj2}, \eqref{3a1}, \eqref{ed1}, and \eqref{ed2}.

Finally,  it follows from \eqref{z1.5}, \eqref{ed2}, and \eqref{ed3} that \eqref{ed4} holds.   The proof of Lemma \ref{lem-ed1} is completed. \hfill $\Box$

 % The proof of Lemma \ref{lem-ed1} is finished. \hfill $\Box$

%\newpage
We will use   Lemmas \ref{lem3.3}--\ref{lem-ed1}   to  prove the following time-independent bound  on the $L^1(0,T;L^\infty)$-norm of $\na u$ which is important for obtaining  the uniform one (with respect to time) on the $L^\infty(0,T;L^q)$-norm of the  gradient of $\mn$.

\begin{lemma}\label{lem5.1} Let $(\n, u, P)$ be a  smooth  solution to  \eqref{1.1}--\eqref{1.3} satisfying \eqref{3a1}.  Then there exists  a generic positive constant  $C$  depending only on  $q,$ $\beta,$ $\ban,$ $\xmu,$ $\smu,$ $\|\n_0\|_{L^{3/2}},$  and $M$ such that
\be\label{3d4.1}    \int_0^T \|\na u\|_{L^\infty}dt \le C  \|u_0\|_{\dot H^{\beta}} . \ee
\end{lemma}

\pf First, it follows from the Gagliardo-Nirenberg inequality that for any $   p\in [2,\min\{6,q\}],$
\bnn\label{nz1.4}%\label{3d1.1}
\ba&
\|\n u_t+\n u\cdot\na u\|_{L^p} \\
   &\le C\|\n^{1/2} u_t\|_{L^2}^{\frac{6-p}{2p}}\|\n^{1/2} u_t\|_{L^6}^{\frac{3p-6}{2p}}+C\|  u\|_{L^6}\|\na u\|_{L^{\frac{6p}{6-p}}}\\
      &\le C\|\n^{1/2} u_t\|_{L^2}^{\frac{6-p}{2p}}\|\na u_t\|_{L^2}^{\frac{3p-6}{2p}}
      + C\|\na u\|_{L^2}\|\na u\|_{L^2}^{\frac{p}{5p-6}}\|\na^2 u\|_{L^p}^{\frac{4p-6}{5p-6}},
      %&\le C(\ban,\|\n_0\|_{L^{3/2}})  \|\n^{1/2} u_t\|_{L^2}+C(\ban,\|\n_0\|_{L^{3/2}}) )\|\na u\|_{L^2}^3+\eta \|\na^2 u\|_{L^2},
\ea\enn
which together with     \eqref{3rd3'}  and  \eqref{z1.5}  gives
\be\ba\label{3d4.2}       \|\na^2u\|_{L^p} +\|\na P\|_{L^p}
 &\le C \|\n u_t+\n u\cdot\na u\|_{L^{6/5}\cap L^p}   \\
&\le   C\|\n^{1/2} u_t\|_{L^2}^{\frac{6-p}{2p}}\|\na u_t\|_{L^2}^{\frac{3p-6}{2p}}+ C\|\na u\|_{L^2}^{\frac{6p-6}{p}}\\
&\quad+\frac12\|\na^2 u\|_{L^p}+C\|\n^{1/2} u_t\|_{L^2}+C \|\na u\|_{L^2}^3 .\ea\ee

Then, setting \be\la{zv1}r\triangleq\frac12\min\xl\{ q+3,\frac{3(5-2\beta)}{3-2\beta}\xr\}\in\xl(3,\min\xl\{ q,\frac{6}{3-2\beta}\xr\}\xr),\ee
one derives  from the Sobolev inequality  and \eqref{3d4.2}    that
\be\label{3d4.6} \ba    \|\na u\|_{L^\infty}  & \le C\|\na  u\|_{L^2}+C\|\na^2 u\|_{L^r}\\ &\le  C\|\na u\|_{L^2}+C \|\n^{1/2}u_t\|_{L^2}+C\|\n^{1/2}u_t\|_{L^2}^\frac{6-r}{2r}\|\na u_t\|_{L^2}^\frac{3r-6}{2r} \\&\quad+C\|\na u\|_{L^2}^\frac{6(r-1)}{r} .\ea\ee

%Hence,
%\be\ba\label{3d4.7}    \int_0^T \|\na u\|_{L^\infty}dt=\int_0^1 \|\na u\|_{L^\infty}dt+\int_1^T \|\na u\|_{L^\infty}dt. \ea\ee

On the one hand, it follows from  \eqref{gj1} and  \eqref{gj2} that for $t\in (0,{\zeta(T)}],$
\bnn\ba  \|\na u \|_{L^\infty}&\le C \|u_0\|_{\dot H^{\beta}}t^\frac{\beta-2}{2} +C \|u_0\|_{\dot H^{\beta}}^\frac{6-r}{2r}t^\frac{\beta-2}{2} \xl(t^{2-\beta}\|\na u_t\|_{L^2}^2\xr)^\frac{3r-6}{4r}\\&\quad +C \|u_0\|_{\dot H^{\beta}}^2t^{2r(\beta-1)/3}+C \|\na u\|_{L^2}^4,\ea\enn
  which together with \eqref{3a1}, \eqref{gj2}, and \eqref{zv1}   gives
\be\ba\label{3d4.8}    & \int_0^{\zeta(T)} \|\na u\|_{L^\infty}dt \\
&\le  C \|u_0\|_{\dot H^{\beta}}+ C \|u_0\|_{\dot H^{\beta}}^\frac{6-r}{2r} \xl(\int_0^1t^\frac{2(\beta-2)r}{r+6}dt\xr)^\frac{r+6}{4r} \xl(\int_0^1t^{2-\beta}\|\na u_t\|_{L^2}^2dt\xr)^\frac{3r-6}{4r}\\
&\le  C \|u_0\|_{\dot H^{\beta}}.\ea\ee

On the other hand, using \eqref{3d4.6}, \eqref{ed2}, and \eqref{ed3}, we obtain that  for $t\in [{\zeta(T)}, T],$
\bnn \ba    \|\na u\|_{L^\infty}    &\le C \|\n^{1/2}u_t\|_{L^2}+C \|\na u_t\|_{L^2} +C\|\na u\|_{L^2} +C\|\na u\|_{L^2}^{6}\\  &\le C\|u_0\|_{\dot H^{\beta}}e^{-\si t/2}  +C \|\na u_t\|_{L^2} ,\ea\enn
%which implies
and thus
\be\ba \label{jia9}
\int_{\zeta(T)}^T\|\na u\|_{L^\infty} dt&\le C\|u_0\|_{\dot H^{\beta}}+C\xl(\int_{\zeta(T)}^Te^{-\si t}dt\xr)^{1/2}\xl(\int_{\zeta(T)}^Te^{ \si t}\|\na u_t\|_{L^2}^2dt\xr)^{1/2}\\&\le C\|u_0\|_{\dot H^{\beta}}.\ea\ee
%due to \eqref{ed3}.
Combining this with   \eqref{3d4.8}  gives \eqref{3d4.1} and   finishes the proof of Lemma \ref{lem5.1}. \hfill$\Box$

With  Lemmas \ref{lem3.3}--\ref{lem5.1}  at hand, we are in a position to prove  Proposition \ref{pr1}.

\emph{Proof of Proposition \ref{pr1}.}  Since $\mn$ satisfies
\bnn\label{3d5.4}   ( \mn)_t+  u\cdot\na  \mn=0,\enn standard calculations show that
\be\label{3d5.5}   \frac{d}{dt}\|\na\mn\|_{L^{q}}\le q\|\na u\|_{L^\infty}\|\na\mn\|_{L^{q}},\ee
which together with Gr\"onwall's inequality and \eqref{3d4.1} yields
\be \la{cbd2z}\ba    \sup_{t\in[0,T]} \|\na\mn\|_{L^{q}}&\le \|\na\mu(\n_0)\|_{L^{q}} \exp\xl\{ q\int_0^T \|\na u\|_{L^\infty}dt \xr\}\\
&\le \|\na\mu(\n_0)\|_{L^{q}} \exp\xl\{  C \|u_0\|_{\dot H^{\beta}}\xr\}\\
&\le 2\|\na\mu(\n_0)\|_{L^{q}},
\ea\ee
provided
\be\ba\label{3d5.6}   \|u_0\|_{\dot H^{\beta}} \le\ve_1\triangleq   C^{-1} \ln 2.
\ea\ee

 Moreover, it follows from \eqref{gj1} and \eqref{ed2} that
\begin{equation}\label{3a1.2} \ba
 \int_0^T\|\na u\|_{L^2}^4dt
 \le&  \sup_{t\in[0,{\zeta(T)}]}\xl(t^{1-\beta}\|\na u\|_{L^2}^2\xr)^2\int_0^{\zeta(T)} t^{2\beta - 2 } dt\\&+\sup_{t\in[{\zeta(T)},T]}\xl(e^{\si t}\|\na u\|_{L^2}^2\xr)^2\int_{\zeta(T)}^Te^{-2\si t}dt\\
 \le & C \|u_0\|_{\dot H^\beta}^4\le \|u_0\|_{\dot H^\beta}^2,
\ea\end{equation}
provided
\be\ba\label{3cb5.6}   \|u_0\|_{\dot H^{\beta}} \le\ve_2\triangleq   C^{-1/2}  .
\ea\ee

Choosing $\ve_0\triangleq \min\{1,\ve_1,\ve_2\},$ we directly obtain \eqref{3a2} from \eqref{cbd2z}--\eqref{3cb5.6}. The proof of Proposition \ref{pr1} is finished.  \hfill$\Box$

The following Lemma \ref{lem5.3} is necessary for further estimates on the higher-order derivatives of the strong solution $(\n,u,P).$
\begin{lemma}\label{lem5.3} Let $(\n, u, P)$ be a  smooth  solution to  \eqref{1.1}--\eqref{1.3} satisfying \eqref{3a1}. Then  there exists a positive constant $C$ depending only on $q,\beta,\ban, \xmu,\smu, M, $ $\|\n_0\|_{L^{3/2}},$      and $\|\na u_0\|_{L^2}$  such that for   $ p_0\triangleq\min\{6,q\} ,$
\be\label{3d6.1}  \ba& \sup_{t\in[0,T]} e^{\sigma t}\xl(\|\na u\|^2_{L^2}+\zeta\|\na u\|_{H^1}^2+\zeta\|P\|_{H^1}^2\xr)+\int_0^T \zeta e^{\sigma t}\|\na u_t\|_{L^2}^2  dt
\\& \quad +\int_0^T e^{\sigma t}\xl( \|\na u\|_{H^1}^2+\|P\|_{H^1}^2 + \zeta\|\na u\|_{  W^{1, p_0}}^2
+ \zeta\|P\|_{  W^{1, p_0}}^2\xr)dt\leq C.\ea\ee
%\be\label{3d7.1}  \ba& \sup_{t\in[0,T]} \|\n\|_{H^2} +\int_0^T\xl(\|\na^3 u\|_{L^2}^2\xr)dt\leq C.\ea\ee
\end{lemma}
\pf First, multiplying  \eqref{3d5} by $ e^{\sigma t},$   we get after using    Gr\"onwall's inequality, \eqref{ed1}, and \eqref{3a2} that
\be\label{3d6.2}   \sup_{t\in[0,T]}   e^{\sigma t}\|\na u\|_{L^2}^2 +\int_0^T   e^{\sigma t}\|\n^{1/2} u_t\|_{L^2}^2 dt\leq C.\ee
Combining this with \eqref{3d2.6} gives
\be\ba\notag\label{3d6.3}&
 \frac{d}{dt}\int\rho|u_t|^{2}dx+\xmu\int|\nabla u_t|^{2}dx   \le C \|\n^{1/2}u_t\|_{L^2}^4+C \|\na u \|_{L^2}^2,
\ea\ee
which along with Gr\"onwall's inequality,    \eqref{3d6.2}, and \eqref{ed1} implies that
\be\label{3d6.4}   \sup_{t\in[0,T]}\zeta e^{\sigma t} \|\n^{1/2} u_t\|_{L^2}^2 +\int_0^T \zeta e^{\sigma t}\|\na u_t\|_{L^2}^2 dt\leq C.\ee
 Combining this,    \eqref{z1.5}, and   \eqref{3d6.2} gives
\be\label{3d6.7}     \ba& \sup_{t\in[0,T]}\zeta e^{\sigma t}\xl(\|\na u\|_{H^1}^2+ \|P\|_{H^1}^2 \xr)+\int_0^Te^{\sigma t}\xl( \|\na u\|_{H^1}^2+ \|P\|_{H^1}^2 \xr)dt \leq C.\ea\ee

Finally,   it follows from  \eqref{z1.5}, \eqref{3d4.2}, \eqref{3d6.2}, and \eqref{gj0} that for $   p_0\triangleq\min\{6,q\} ,$
\be\ba\label{3d6.5}     &\|\na u\|_{H^1\cap W^{1,p_0}} +\|P\|_{H^1\cap W^{1,p_0}}
%&= \|\na u\|_{L^p} +\|\na^2 u\|_{L^p} +\|P\|_{L^p} +\|\na P\|_{L^p}
  \le  C\|\na u_t\|_{L^2}  +C \|\na u \|_{L^2},\ea\ee
which together with  \eqref{3d6.4} and  \eqref{ed1}  implies
\bnn\ba &\int_0^T\zeta e^{\sigma t}\xl(\|\na u\|^2_{  W^{1,p_0}} + \|P\|^2_{ W^{1,p_0}} \xr)dt \le C .\ea\enn
This combined with \eqref{3d6.2}--\eqref{3d6.7} gives \eqref{3d6.1} and completes the proof of Lemma \ref{lem5.3}. \hfill$\Box$

The following Proposition \ref{lem5.a3} is concerned with the  estimates on the higher-order derivatives of the strong solution $(\n,u,P)$ which in particular imply the continuity in time of both $\na^2 u$ and $\na P$ in the $L^2\cap L^p$-norm.

\begin{proposition}\label{lem5.a3} Let $(\n, u, P)$ be a  smooth  solution to  \eqref{1.1}--\eqref{1.3} satisfying \eqref{3a1}. Then  there exists a positive constant $C$ depending only on $q,\beta,\ban, \xmu,\smu,  \|\n_0\|_{L^{3/2}},  M, $ $\|\na u_0\|_{L^2}$,    and $\|\na\n_0\|_{L^2} $  such that for $ p_0\triangleq\min\{6,q\} $ and $ q_0\triangleq 4q/(q-3)$,
\be\label{3dr6.1}  \ba& \sup_{t\in[0,T]}\zeta^{q_0 }   e^{\sigma t}\xl(\|\na u\|^2_{W^{1,p_0}}+ \|P\|^2_{W^{1,p_0}}+ \|\na u_t\|^2_{L^2}\xr)   \\&
 +\int_0^T\zeta^{q_0+1} e^{\sigma t}\xl(  \|(\n  u_{t})_t \|_{L^2}^2+  \|\na u_t\|_{  L^{ p_0}}^2
+  \|P_t\|_{L^2\cap L^{ p_0}}^2\xr)dt\leq C.\ea\ee %where $ p_0\triangleq\min\{6,q\} $ and $ q_0\triangleq 4q/(q-3).$
\end{proposition}

{\it Proof. } First,   similar to \eqref{3d5.5} and \eqref{cbd2z}, we have
\be\ba\label{3d6.21}   \sup_{0\le t\le T}\|\na\n\|_{L^2 } \le 2\|\na\n_0\|_{L^2 },
\ea\ee   which together with the Sobolev inequality and \eqref{3d6.2} gives
\be \la{ij6}\ba\|\n_t\|_{L^2\cap L^{3/2}}&=\|u\cdot\na\n\|_{L^2\cap L^{3/2}}\\&\le C\|\na\n\|_{L^2}\|\na u\|_{L^2}^{1/2}\|\na u\|_{H^1}^{1/2} \le C\|\na u\|_{H^1}^{1/2} .\ea\ee

Next, it follows from \eqref{utt} that $u_t$ satisfies
\bnn\label{5rd3}
\begin{cases}
 -\div(2\mn d_t)  +\nabla P_t=\ti F+\div g,\\
 \div u_t=0,
\end{cases}
\enn  with \bnn  \ti F\triangleq  -
\n u_{tt}- \n u \cdot\na u_t-\n_tu_t-(\n u)_t\cdot\na u, \quad g\triangleq-2u\cdot\na\mn  d . \enn
 Hence, one can deduce from   Lemma \ref{stokes} and the Sobolev  inequality that
 \be \la{3'd.1}\ba  \|\na u_t\|_{ L^2\cap L^{p_0}}+\|P_t\|_{ L^2\cap L^{p_0}}     \le C\|\ti F\|_{L^{6/5}\cap L^{\frac{3p_0}{p_0+3}}}+C\|g\|_{L^2\cap L^{p_0}}.\ea\ee
Using \eqref{3a1}, \eqref{gj0}, \eqref{ij6},  \eqref{3d6.1},  and \eqref{3d6.5}, we get by direct calculations   that \be\la{ij1}\ba   &\|\ti F\|_{L^{6/5}\cap L^{\frac{3p_0}{p_0+3}}}  \\ &\le  C\|\n\|^{1/2}_{L^{3/2}\cap L^{\frac{3p_0}{6-p_0}}}\|\n^{1/2} u_{tt} \|_{L^2}+C\|\n\|_{L^{3}\cap L^{\frac{6p_0}{6-p_0}}}\|u\|_{L^\infty}\|\na u_t\|_{L^2}\\&\quad+ C\| \n_t\|_{L^2\cap L^{3/2}} \xl(\|u_t\|_{L^6\cap L^{\frac{6p_0}{6-p_0}} } +\|\na u\|_{H^1}^2+\|\na u\|_{H^1} \|\na u\|_{W^{1,p_0}}\xr) \\&\quad+C\|\n\|_{L^2\cap L^{p_0}}\| u_t\|_{L^6}\|\na u \|_{L^6}\\ &\le
 C \| \sqrt{\n}  u_{tt} \|_{L^2}+\ve \|\na u_t\|_{L^{p_0}}+C (\ve)\|\na u_t\|_{L^2}(1+\|\na u\|^{3/2}_{H^1} ) +C\|\na u\|_{H^1}^{5/2},\ea\ee and that
\be\la{ij2}\ba \|g\|_{L^2\cap L^{p_0}}&\le C\|\na \mn\|_{L^q}\| u\|_{L^6\cap L^\infty}\|\na u\|_{L^2\cap L^\infty}\\& \le  C\|\na u_t\|_{L^2} \|\na u\|_{H^1}+C \|\na u\|_{H^1}^2,\ea\ee
where in the second inequality one has used the following simple fact
\be\la{ij3}\|\na u\|_{L^\infty}\le C\|\na u\|_{H^1\cap W^{1,p_0}}\le C\|\na u_t\|_{L^2}+C\|\na u\|_{L^2},\ee due to the Sobolev  inequality and \eqref{3d6.5}.
Then,  putting  \eqref{ij1} and \eqref{ij2} into \eqref{3'd.1}, we obtain after choosing $\ve$ suitably small that
 \be \la{3d.1}\ba &\|\na u_t\|_{ L^2\cap L^{p_0}}+\|P_t\|_{ L^2\cap L^{p_0}}  \\ &\le
 C \| \sqrt{\n}  u_{tt} \|_{L^2} +C \|\na u_t\|_{L^2}(1+\|\na u\|^2_{H^1} ) +C\|\na u\|_{H^1}+C\|\na u\|_{H^1}^3.\ea\ee

Now, multiplying \eqref{utt} by  $u_{tt} $ and integrating the resulting equality by parts lead to
\be\la{utt1}\ba & \int \n |u_{tt} |^2dx+ \frac{d}{dt}\int \mu(\n)|d_t|^2dx\\&= \int \div(\mu(\n)u)|d_t|^2dx -\int \n( u \cdot \na u_t+ u_t\cdot\na u)\cdot u_{tt} dx  -\int \n_t u^j_t u^j_{tt}dx\\&\quad -\int \n_t u \cdot\na u^j   u^j_{tt}dx  -2\int\pa_i(u^k\pa_k\mu(\n)d^j_i)u^j_{tt}dx  \triangleq \sum_{i=1}^{5}I_i.\ea\ee

We will use \eqref{3d6.1}, \eqref{3d.1},  and the Sobolev inequality to estimate each term on the righthand side of \eqref{utt1} as follows:

First, it follows from \eqref{3a1}, \eqref{3d6.1}, and  \eqref{3d.1} that
\be\la{ij4}\ba|I_1|&\le C\|u\|_{L^\infty}\|\na \mn\|_{L^q}\|\na u_t\|^{\frac{2(p_0q-p_0-2q)}{q(p_0-2)}}_{L^2}\|\na u_t\|^{\frac{2 p_0 }{q(p_0-2)}}_{L^{p_0}}\\&\le \ve\|\na u_t\|^2_{L^{p_0}}+C(\ve)\|\na u\|_{H^1}^{\frac{q(p_0-2)}{ p_0q-p_0-2q }} \|\na u_t\|^2_{L^2} \\&\le C\ve \|\sqrt{\n} u_{tt}\|_{L^2}^2+C(\ve)(1+\|\na u\|_{H^1}^{ q_0})\|\na u_t\|^2_{L^2}  \\&\quad + C(\ve)\| \na  u\|_{H^1}^2+ C(\ve)\| \na  u\|_{H^1}^6, \ea\ee where in the last inequality we have used $$\frac{q(p_0-2)}{ p_0q-p_0-2q }\in [1, q_0].$$

Next, H\"older's inequality gives
\be\ba |I_2|\le \ve \int \n |u_{tt} |^2dx+ C(\ve)\| \na  u\|_{H^1}^2\|\na u_t\|_{L^2}^2.\ea\ee

Then, direct calculations show
\be\la{ij5}\ba I_3&=-\frac12 \frac{d}{dt}\int \n_t|u_t|^2dx+\int (\n u^i)_t \pa_iu^j_tu^j_tdx \\ &\le-\frac12 \frac{d}{dt}\int \n_t|u_t|^2dx+C\|\n\|_{L^6}\|\na u_t\|_{L^2}\|u_t\|_{L^6}^2\\&\quad + C\|\n_t\|_{L^2}\|u \|_{L^\infty}  \|\na u_t\|_{L^3}\|u_t\|_{L^6}  \\ &\le -  \frac{d}{dt}\int \n u\cdot\na  u_t^j  u_t^jdx+C(\ve)(1+\|\na u_t\|_{L^2} + \|\na u \|^4_{H^1} )\|\na u_t\|_{L^2}^2\\&\quad +\ve \int\n|u_{tt}|^2dx+C(\ve)\|\na u \|^2_{H^1}+C(\ve)\|\na u \|^6_{H^1}, \ea\ee
where in the last inequality one has used \eqref{ij6} and \eqref{3d.1}.

Next, it follows from \eqref{1.1}$_1$ and  \eqref{ij6}  that
\be\ba I_4=&  -\frac{d}{dt}\int \n_t u \cdot\na u^j   u^j_{t }dx+\int (\n u^ i)_{t }\pa_i( u \cdot\na u^j   u^j_{t })dx+\int \n_{t } (u \cdot\na u^j )_t  u^j_{t }dx\\=&-\frac{d}{dt}\int \n_t u \cdot\na u^j   u^j_{t }dx+\int \n u^i_{t }( u \cdot\na u^j  \pa_i u^j_{t }+\pa_i( u \cdot\na u^j )  u^j_{t })dx\\&+\int \n_{t }  u^ i ( u \cdot\na u^j  \pa_i u^j_{t }+\pa_i( u \cdot\na u^j )  u^j_{t }) dx  +\int \n_{t }(u \cdot\na u^j )_t  u^j_{t } dx\\\le &-\frac{d}{dt}\int \n_t u \cdot\na u^j   u^j_{t }dx+C\|u_t\|_{L^6}\|\na u\|_{H^1}^2(\|\na u_t\|_{L^2}+\|u_t\|_{L^6})\\&+C\|\n_t\|_{L^2}\|\na u\|^{1/2}_{H^1} (\|\na u_t\|_{L^2}\|\na u\|_{H^1} \|\na u\|_{H^1\cap W^{1,p_0}}+\|u_t\|_{L^6}\|\na u\|_{H^1}^2)\\&+C\|\n_t\|_{L^2}\|u_t\|_{L^6}(\|u_t\|_{L^6}\|\na u\|_{H^1} +\|\na u_t\|_{L^3}\|\na u\|_{H^1} )\\\le &-\frac{d}{dt}\int \n_t u \cdot\na u^j   u^j_{t }dx+C(\ve)(1+\|\na u_t\|_{L^2} + \|\na u \|^4_{H^1} )\|\na u_t\|_{L^2}^2\\&\quad +\ve \int\n|u_{tt}|^2dx+C(\ve)\|\na u \|^2_{H^1}+C(\ve)\|\na u \|^6_{H^1}. \ea\ee

Finally, direct calculations lead to
\be\la{ij51}\ba  I_5&= -2\frac{d}{dt} \int\pa_i(u^k\pa_k\mu(\n)d^j_i)u^j_{t }dx- 2\int\pa_i( u^k\mu (\n) \pa_k d^j_i )_tu^j_{t}dx\\&\quad+2\int\pa_i(u^k\pa_k(\mu(\n)d^j_i) )_tu^j_{t}dx\\& =2 \frac{d}{dt} \int u^k\pa_k\mu(\n)d^j_i \pa_iu^j_{t }dx+2\int(\mu (\n)  u^k\pa_k d^j_i)_t \pa_iu^j_{t }dx\\&\quad-2\int(\pa_iu^k \mu(\n)d^j_i  )_t\pa_ku^j_{t}dx-2\int u^k_t\pa_{i}(\mu(\n)d^j_i)  \pa_ku^j_{t}dx \\&\quad-2\int u^k(\pa_{i}(\mu(\n)d^j_i) )_t\pa_ku^j_{t}dx \\&\triangleq 2 \frac{d}{dt} \int u^k\pa_k\mu(\n)d^j_i \pa_iu^j_{t }dx+ \sum_{i=1}^{4}I_{5,i}.\ea\ee
We estimate each $I_{5,i} (i=1,\cdots,4)$  as follows:

First, integration by parts gives
 \be\ba I_{5,1}&=2\int(\mu (\n)  u^k)_t\pa_k  d^j_i\pa_iu^j_{t }dx+2\int \mu (\n)  u^k\pa_k  (d^j_i)_t\pa_iu^j_{t }dx\\&=-2\int u\cdot\na\mu (\n)   u^k\pa_k  d^j_i\pa_iu^j_{t }dx+2\int \mu (\n)  u^k_t\pa_k  d^j_i\pa_iu^j_{t }dx\\&\quad- \int \div(\mu (\n)  u) |d_t|^2dx\\&\le C \|u\|_{L^{6q/(q-3)}}^2\|\na\mn\|_{L^q} \|\na^2 u\|_{L^3}\|\na u_t\|_{L^3}\\&\quad+ C   \|\na^2 u\|_{L^3} \|\na u_t\|_{L^2}^2+|I_1|  \\&\le C\ve \|\n^{1/2}u_{tt}\|_{L^2}^2+C(\ve)(1+\|\na u\|_{H^1}^{q_0}+\|\na u_t\|_{L^2})\|\na u_t\|^2_{L^2} \\&\quad+ C(\ve)\| \na  u\|_{H^1}^6+ C(\ve)\| \na  u\|_{H^1}^2 ,\ea\ee
where in the last inequality we have used \eqref{3d6.1}, \eqref{3d6.5}, \eqref{3d.1}, and \eqref{ij4}.

Then, it follows from \eqref{3a1} and \eqref{ij3} that
\be\ba |I_{5,2}| &\le C \|u\|_{L^\infty}\| \na\mu(\n) \|_{L^q}\|\na u\|_{L^{3q/(q-3)}}\|\na u\|_{L^6}\|\na u_t\|_{L^2}\\&\quad+C\|\na u\|_{L^\infty} \|\na u_t\|_{L^2}^2 \\ &\le C   \|\na u\|_{H^1\cap W^{1,p_0}}\xl(\|\na u\|_{H^1}^2\|\na u_t\|_{L^2} +C \|\na u_t\|_{L^2}^2 \xr) \\&\le C  \|\na u\|_{H^1}^4 +C(1+\|\na u\|_{H^1}^2+\|\na u_t\|_{L^2})\|\na u_t\|^2_{L^2} .\ea\ee
Similarly, combining H\"older's inequality and \eqref{3d6.5} leads to
\be\ba |I_{5,3}|&\le C\|u_t\|_{L^6}\|\na u_t\|_{L^2}(\|\na\mn\|_{L^q}\|\na u\|_{L^{3q/(q-3)}}+\|\na^2 u\|_{L^3})\\ &\le C \|\na u_t\|_{L^2}^2(\|\na u_t\|_{L^2}+\|\na u\|_{H^1}) .\ea\ee

Finally, using \eqref{1.1}$_2$ and \eqref{1.1}$_3$, we obtain after integrating by parts that
\be
\la{ij16}\ba I_{5,4}&= -2\int u^k  \pa_j P_t\pa_{ k}u^j_{t}dx  -2
  \int u^k(\n u^j_t+\n u\cdot\na u^j )_t\pa_{ k}u^j_{t}dx\\&= 2\int \pa_ju^k   P_t\pa_{ k}u^j_{t}dx -2
  \int u^k \n u^j_{tt} \pa_{ k}u^j_{t}dx \\&\quad -2
  \int u^k(\n_t u^j_t+(\n u\cdot\na u^j )_t)\pa_{ k}u^j_{t}dx  \\&\le C\|\na u \|_{L^6}\|P_t\|_{L^3} \|\na u_t\|_{L^2} +C\|\sqrt{\n}  u_{tt} \|_{L^2} \|\na u\|_{H^1}  \|\na u_t\|_{L^2} \\& \quad+C\|u \|_{L^\infty}\|\na u_t\|_{L^2} \|\n_t\|_{L^2}( \|u_t\|_{L^\infty}+\|\na u\|_{H^1} \|\na u\|_{L^\infty} )  \\& \quad +C\|u \|_{L^\infty}\|\na u_t\|_{L^2}( \|\na u_t\|_{L^2} +\|
    u_t\|_{L^6} )\|\na u\|_{H^1}   \\&\le C\ve \int \n |u_{tt} |^2dx+ C(\ve) (1+\| \na u_t\|_{L^2}+\| \na u\|_{H^1}^4) \|\na u_t\|_{L^2}^2 \\&\quad +C(\ve)\|\na u\|_{H^1}^2+C(\ve)\|\na u\|_{H^1}^6 , \ea\ee where in the last inequality one has used \eqref{3d.1} and \eqref{ij6}.

Substituting \eqref{ij4}--\eqref{ij16} into \eqref{utt1}, we get after choosing $\ve$ suitably small that

\be\la{ij20}\ba & \frac{d}{dt}\int \mu(\n)|d_t|^2dx+\Psi' (t) +\frac12\int \n |u_{tt}|^2dx \\ &\le C  (1+\| \na u_t\|_{L^2}+\| \na u\|_{H^1}^{q_0}) \|\na u_t\|_{L^2}^2    +C \|\na u\|_{H^1}^2+C \|\na u\|_{H^1}^6,\ea\ee
where
\bnn \la{ij21}\ba  \Psi(t)\triangleq    -\int \n u\cdot\na  u_t^j  u_t^jdx-\int \n_t u \cdot\na u^j   u^j_{t }dx +2\int u^k\pa_k\mu(\n)d^j_i \pa_iu^j_{t }dx\ea\enn satisfies \be\la{ij22} \ba |\Psi(t)|\le& C\|\sqrt{\n}u_t\|_{L^2}\|\na u_t\|_{L^2}\|\na u\|_{H^1}+C\|\n_t\|_{L^2}\| u\|_{L^6}\|u_t\|_{L^6}\|\na u\|_{L^6}\\&+C \|\na\mn\|_{L^q}\|\na u_t\|_{L^2}\|\na u\|_{H^1}^2\\ \le& \frac{1}{4}\xmu\|\na u_t\|_{L^2}^2+C\|\sqrt{\n}u_t\|_{L^2}^2 \|\na u\|_{H^1}^2+C\|\na u\|_{H^1}^4,\ea\ee due to \eqref{3a1} and \eqref{ij6}.

Then,  multiplying  \eqref{ij20} by $\zeta^{q_0} e^{\sigma t} $  and
noticing that \eqref{3d6.1} gives
\bnn \la{ij23}\ba \zeta^{q_0}(1+\| \na u_t\|_{L^2}+\| \na u\|_{H^1}^{q_0}) \|\na u_t\|_{L^2}^2 \le C  \zeta^{q_0+1} \|\na u_t\|_{L^2}^4+C\zeta \| \na u_t\|^2_{L^2},\ea\enn   we get after using    Gr\"onwall's inequality, \eqref{ij22}, \eqref{3d6.1}, and \eqref{3d6.4} that
\be \la{ij24}\sup_{0\le t\le T}\zeta^{q_0}e^{\sigma t}\|\na u_t\|_{L^2}^2+\int_0^T \zeta^{q_0} e^{\sigma t}\int\n |u_{tt}|^2dxdt\le C.\ee
Furthermore, it follows from   \eqref{ij6} and \eqref{3d6.1} that
 \bnn \ba\|(\n u_t)_t\|_{L^2}^2\le C\|\na u\|_{H^1}\|\na u_t\|_{L^2
\cap L^{p_0}}^2+C\|\n^{1/2}u_{tt}\|_{L^2}^2,\ea \enn
which together with \eqref{ij24}, \eqref{3d6.5}, \eqref{3d.1},  and \eqref{3d6.1} gives \eqref{3dr6.1} and thus completes the proof of Proposition \ref{lem5.a3}.   \hfill $\Box$

\section{Proofs of Theorems \ref{thm1} and \ref{thm2}}

With all the a priori estimates in Section 3 at hand, we are now in a position to  prove Theorems \ref{thm1} and \ref{thm2}.

\emph{Proof of Theorem \ref{thm1}.} First,
by Lemma 2.1, there exists a $T_{*}>0$ such that the Cauchy problem \eqref{1.1}-\eqref{1.3} has a unique local strong solution $(\rho,u,P)$ on $\mathbb{R}^3\times(0,T_{*}]$. It follows from \eqref{2.2}  that there exists a $T_1\in (0, T_*]$ such that \eqref{3a1} holds for $T=T_1$.

Next, set
\begin{equation}\label{20.1}
T^{*}\triangleq\sup \{T | (\n, u, P) \mbox{ is a strong solution on } \rrr\times(0,T] \mbox{ and }\eqref{3a1}\ \text{holds}\}.
\end{equation}
Then $T^*\ge T_1>0$. Hence,  for any $0< \tau<T\leq T^{*}$ with $T$ finite, one deduces from \eqref{3d6.1} and \eqref{3dr6.1}  that
\begin{equation}\label{20.2}
 \nabla u, \  P\in C([\tau,T];L^2)\cap C(\overline \rrr\times [\tau,T] ) ,
\end{equation}
where one has used the standard embedding
\begin{equation*}
L^{\infty}(\tau,T;H^1\cap W^{1,p_0})\cap H^1(\tau,T;L^2)\hookrightarrow C([\tau,T];L^2)\cap C(\overline \rrr\times [\tau,T] ) .
\end{equation*}
Moreover, it follows from \eqref{3a1}, \eqref{gj0}, \eqref{3d6.21}, and \cite[Lemma 2.3]{L1996} that
\begin{equation}\label{20.3}
\rho\in C([0,T];L^{3/2}\cap H^1), \quad\na\mn \in C([0,T];L^q).
\end{equation}
Thanks to \eqref{3d6.2} and \eqref{3dr6.1},
the standard arguments yield that
\bnn\label{20.4}
\rho   u_t \in H^1(\tau,T;L^2)\hookrightarrow C([\tau,T];L^2),
\enn
which together with \eqref{20.2} and \eqref{20.3} gives
\begin{equation}\label{20a4}
\rho   u_t+\n u\cdot\na u \in  C([\tau,T];L^2).
\end{equation} Since $(\n,u)$ satisfies \eqref{3rd7} with $F\equiv \n u_t+\n u\cdot\na u$, we deduce from  \eqref{1.1}, \eqref{20.2},  \eqref{20.3},   \eqref{20a4},  and \eqref{3dr6.1} that
\begin{equation}\label{n20.4}
 \na u,~P\in C([\tau,T];D^1\cap D^{1,p}),
\end{equation} for any $p\in [2,p_0).$

Now, we claim that
\begin{equation}\label{20.5}
T^*=\infty.\end{equation}
Otherwise, $T^*<\infty$.  Proposition \ref{pr1} implies that \eqref{3a2} holds at $T=T^*$. It follows from \eqref{20.2}, \eqref{20.3},  and \eqref{n20.4}  that
$$(\n^*,u^*)(x)\triangleq(\n, u)(x,T^*)=\lim_{t\rightarrow T^*}(\n, u)(x,t)$$
satisfies \bnn \n^*\in L^{3/2}\cap H^1,\quad u^*\in D^{1}_{0,\sigma}\cap D^{1,p}\enn  for any $p\in [2,p_0).$   Therefore, one can  take $(\n^*,\n^*u^*)$ as the initial data and apply
Lemma 2.1 to extend the local strong solution beyond $T^*$.
 This contradicts the assumption of $T^*$ in \eqref{20.1}.
 Hence, \eqref{20.5} holds. We thus finish
the proof of Theorem \ref{thm1} since \eqref{oiq1} and   \eqref{e} follow directly from \eqref{3d6.21} and  \eqref{3dr6.1}, respectively.
 \hfill $\Box$

\emph{Proof of Theorem \ref{thm2}}: With the global existence result at hand (see Proposition \ref{thm3}), one can modify slightly the proofs of Lemma \ref{lem-ed1} and \eqref{3d6.21} to obtain    \eqref{oiq2} and \eqref{eq}.
 \hfill$\Box$

%\emph{Proof of Theorem \ref{thm2}}:  With the global existence result in Theorem \ref{thm3} at hand, one can reprove Lemmas \ref{lem3.3} and \ref{lem3.4}--\ref{lem-ed3} by  the same arguments as those in Section \ref{sec3} step by step, and thus deduce  the large time exponential decay estimates \eqref{hwe2}. \hfill$\Box$

%\newpage


\begin{thebibliography}{10}
\bibitem{zhang1}
H. Abidi, G. L. Gui, P. Zhang, On the decay and stability to global solutions of the 3-D inhomogeneous Navier-Stokes equations, \emph{Comm. Pure Appl. Math.}, \textbf{64}(2011), 832--881.


%\bibitem{zhang4}
%H. Abidi, G.L. Gui, P. Zhang, On  the well-posedness of  three  dimensional  inhomogeneous Navier-Stokes equations  in the critical spaces, \emph{Arch. Ration. Mech. Anal.}, \textbf{204}(2012), 189--230.
\bibitem{zhang5}
H. Abidi, P. Zhang, On the global well-posedness of 2-D density-dependent Navier-Stokes system with variable
viscosity, \emph{J. Differential Equations}, {\bf 259} (2015),  3755--3802.

\bibitem{zhang6}
H. Abidi, P. Zhang, Global well-posedness of 3-D density-dependent Navier-Stokes system with variable viscosity,
\emph{Sci. China Math.}, \textbf{58(6)}(2015), 1129--1150.


\bibitem{AK1973}
S. A. Antontesv,  A. V. Kazhikov,    Mathematical study of flows of nonhomogeneous fluids, Lecture notes (1973), Novosibirsk State University, Novosibirsk, U.S.S.R.



\bibitem{stein}
J. Bergh, J. Lofstrom,  Interpolation spaces, An introduction. Springer-Verlag, Berlin-Heidelberg-New York, 1976.

\bibitem{chen1} Z. M. Chen, A sharp decay result on strong solutions of the Navier-
Stokes equations in the whole space, \emph{ Comm. Partial Differential Equations},
{\bf 16} (1991),   801--820.

\bibitem{kim2004}
Y. Cho, H. Kim, Unique solvability for the density-dependent Navier-Stokes equations, \emph{Nonlinear Anal.}, \textbf{59}(4)(2004),  465--489.


\bibitem{kim2003}
H. Y. Choe, H. Kim, Strong solutions of the Navier-Stokes equations for nonhomogeneous incompressible fluids,
\emph{Comm. Partial Differential Equations}, \textbf{28}(2003), 1183--1201.


\bibitem{huang13}
W. Craig, X. D. Huang, Y. Wang, Global strong solutions for 3D nonhomogeneous incompressible Navier-Stokes equations, \emph{J. Math. Fluid Mech.}, \textbf{15}(2013), 747--758.


%\bibitem{dan1}
%R. Danchin, Density-dependent incompressible viscous fluids in critical spaces, \emph{Proc. Roy. Soc. Edinburgh Sect. A}, \textbf{133}(2003),  1311--1334.



\bibitem{dan2}
R. Danchin, Local and global well-posedness results for flows of inhomogeneous viscous fluids, \emph{Adv. Differential Equations}, \textbf{9}(2004), 353--386.



\bibitem{dan3}
R. Danchin, P. B. Mucha, A Lagrangian approach for the incompressible Navier-Stokes equations with variable density, \emph{Comm. Pure Appl. Math.}, \textbf{65}(2012), 1458--1480.

%\bibitem{dan4}
%R. Danchin, P.B. Mucha, Incompressible flows with piecewise constant density, \emph{Arch. Ration. Mech. Anal.}, \textbf{207(3)}(2013),  991--1023.

\bibitem{D1997}
B. Desjardins, Regularity results for two-dimensional flows of multiphase viscous fluids, \emph{Arch. Rational Mech. Anal.}, \textbf{137}(1997), 135--158.

\bibitem{ft1964} H. Fujita,   T. Kato, On the Navier-Stokes initial value problem. I. \emph{Arch. Rational Mech. Anal.}
{\bf 16}(1964), 269--315.
%\bibitem{H2014} C. He,  X. Huang, Y. Wang,  On some new global existence results for 3D magnetohydrodynamic equations. \emph{Nonlinearity}, \textbf{27} (2014), no. 2, 343--352.

\bibitem{N1959} G. P. Galdi,  An introduction to the mathematical theory of the Navier-Stokes equations. Steady-state problems. Second edition.   Springer, New York, 2011.

\bibitem{H2002}   C. He, L. Hsiao,  The decay rates of strong solutions for Navier-Stokes equations, \emph{J. Math. Anal. Appl.},  \textbf{268}(2002),  417--425.


\bibitem{Hoh4}
D. Hoff,
   Compressible flow in a half-space with Navier boundary
  conditions. \textit{J. Math. Fluid Mech.},  \textbf{7}(2005), no. 3, 315--338.




\bibitem{zhang3}
J. C. Huang, M. Paicu, P. Zhang, Global solutions to 2-D inhomogeneous Navier-Stokes system with general velocity, \emph{J. Math. Pures Appl.},  \textbf{100}(2013), 806--831.


\bibitem{HLX2012}
X. D. Huang,  J. Li,  Z. P. Xin,  Global well-posedness of classical solutions with large oscillations and vacuum to the three-dimensional isentropic compressible Navier-Stokes equations, \emph{Comm. Pure Appl. Math.},  \textbf{65}(2012), 549--585.


%\bibitem{hwjde1}
%X.D. Huang, Y. Wang, Global strong solution to the 2D nonhomogeneous incompressible MHD system, \emph{J. Differential Equations}, \textbf{254}(2014), 511--527.


\bibitem{HW2014}
 X. D. Huang,   Y. Wang,  Global strong solution with vacuum to the two-dimensional density-dependent Navier-Stokes system, \emph{SIAM J. Math. Appl.}, \textbf{46}(2014), 1771--1788.




\bibitem{huang15}
X. D. Huang, Y. Wang, Global strong solution of 3D inhomogeneous Navier-Stokes equations with density-dependent viscosity, \emph{J. Differential Equations}, \textbf{259}(2015), 1606--1627.


\bibitem{Kt1984} T. Kato, Strong $L^p$-solutions of the Navier-Stokes equation in $R^m$, with applications to weak
solutions. \emph{Math. Z.}, {\bf 187}(1984), no. 4, 471--480.

\bibitem{K1974}
A. V. Kazhikov,    Resolution of boundary value problems for nonhomogeneous viscous fluids, \emph{Dokl. Akad. Nauk.}, \textbf{216}(1974), 1008--1010.

%\bibitem{Kt01}H. Koch,  D.  Tataru, Well-posedness for the Navier-Stokes equations. \emph{ Adv. Math.} {\bf 157} (2001), no. 1, 22--35.
%\bibitem{decay16} D.Q. Khai, N.M. Tri, The existence and space-time decay rates of strong solutions to
%Navier-Stokes Equations in weighed $L^\infty(|x|^\gamma dx)\cap L^\infty(|x|^\beta dx)$ spaces, http://arxiv.org/pdf/1601.01723v1.


\bibitem{lady}
O. Ladyzhenskaya, V.A. Solonnikov, Unique solvability of an initial and boundary value problem for viscous incompressible nonhomogeneous fluids, \emph{J. Soviet Math.}, \textbf{9}(1978), 697--749.


%\bibitem{LL2014} J. Li,   Z.L. Liang, On local classical solutions to the Cauchy problem of the two-dimensional barotropic compressible Navier-Stokes equations with vacuum, \emph{J. Math. Pures Appl.}, \textbf{102}(2014),  640--671.


%\bibitem{LX2014} J. Li,  Z.P. Xin,    Global well-posedness and large time asymptotic behavior of classical solutions to the compressible Navier-Stokes equations with vacuum,  http://arxiv.org/abs/1310.1673, 2013.


%\bibitem{liang15}
%Z.L. Liang, Local strong solution and blow-up criterion for the 2D nonhomogeneous incompressible fluids, \emph{J. Differential Equations}, \textbf{7}(2015), 2633--2654.


\bibitem{L1996}
P. L. Lions, Mathematical topics in fluid mechanics, vol. I: incompressible models,  Oxford University Press, Oxford, 1996.

\bibitem{lvshizh1}
B.Q. L\"u, Z.H. Xu, X. Zhong,
On local strong solutions to the Cauchy problem of two-dimensional density-dependent Magnetohydrodynamic equations with vacuum,
 http://arxiv.org/abs/1506.02156

\bibitem{lvshizh}
B.Q. L\"u, X. D. Shi, X. Zhong, Global existence and large time asymptotic behavior of strong solutions to the Cauchy problem of 2D density-dependent Navier-Stokes equations with vacuum, http://arxiv.org/abs/1506.03143.



\bibitem{hls}
  B.Q. L\"u, S.S. Song,
On local strong solutions to the three-dimensional nonhomogeneous incompressible Navier-Stokes equations with density-dependent viscosity and vacuum, submitted.
%\bibitem{prst1} G. Ponce,  R. Racke, T. C. Sideris,  E. S.  Titi,
%Global stability of large solutions to the 3D Navier-Stokes equations.
%Comm. Math. Phys. \textbf{159} (1994), no. 2, 329--341.


\bibitem{sch}M. E. Schonbek, Large time behaviour of solutions to the Navier-
Stokes equations in $H^m$ spaces,\emph{ Comm. Partial Differential Equations},
{\bf 20}(1995),   103--117.
%\bibitem{zhang2}
%M. Paicu, P. Zhang, Global solutions to the 3-D incompressible inhomogeneous Navier-Stokes system, \emph{J. Funct. Anal.}, \textbf{262}(2012),  3556--3584.
%\bibitem{pzz} Paicu, Marius; Zhang, Ping; Zhang, Zhifei Global unique solvability of inhomogeneous Navier-Stokes equations with bounded density. Comm. Partial Differential Equations 38 (2013), no. 7, 1208--1234.
\bibitem{S1990}
 J. Simon,
 Nonhomogeneous viscous incompressible fluids: existence of velocity, density, and pressure,
 \emph{SIAM J. Math. Anal.}, \textbf{21}(1990), 1093--1117.


\bibitem{zjw}
J. W. Zhang, Global well-posedness for the incompressible Navier-Stokes equations with density-dependent viscosity coefficient,  \emph{J. Differential Equations}, \textbf{259}(2015), 1722--1742.


\end{thebibliography}
\end{document}